\documentclass{article}

\newcommand{\mK}{\mathbb{K}}

\usepackage{proof}
\usepackage{latexsym}
\usepackage{amsfonts}
\usepackage{amssymb}
\usepackage{cite}
\usepackage{color}

\newcommand{\brem}{\begin{remark}}
\newcommand{\erem}{\end{remark}}
\newcommand{\blem}{\begin{lemma}}
\newcommand{\elem}{\end{lemma}}
\newcommand{\bth}{\begin{theorem}}
\newcommand{\ethm}{\end{theorem}}
\newcommand{\benu}{\begin{enumerate}}
\newcommand{\eenu}{\end{enumerate}}
\newcommand{\bdes}{\begin{description}}
\newcommand{\edes}{\end{description}}
\newcommand{\bdf}{\begin{definition}}
\newcommand{\edf}{\end{definition}}
\newcommand{\bcor}{\begin{cor}}
\newcommand{\ecor}{\end{cor}}
\newcommand{\bprp}{\begin{proposition}}
\newcommand{\eprp}{\end{proposition}}
\newcommand{\bmlem}{\begin{mlemma}}
\newcommand{\emlem}{\end{mlemma}}
\newcommand{\bclm}{\begin{claim}}
\newcommand{\eclm}{\end{claim}}
\newcommand{\bprf}{{\bf Proof}.\hspace{2mm}}

\newcommand{\eprf}{\hspace*{\fill} $\Box$}

\newcommand{\beqn}{\begin{equation}}
\newcommand{\eeqn}{\end{equation}}
\newcommand{\beqnarr}{\begin{eqnarray}}
\newcommand{\eeqnarr}{\end{eqnarray}}
\newcommand{\beqnarrs}{\begin{eqnarray*}}
\newcommand{\eeqnarrs}{\end{eqnarray*}}

\newcommand{\spand}{\,\&\,}

\newtheorem{theorem}{Theorem}[section]
\newtheorem{definition}[theorem]{Definition}
\newtheorem{proposition}[theorem]{Proposition}
\newtheorem{lemma}[theorem]{Lemma}
\newtheorem{cor}[theorem]{Corollary}
\newtheorem{remark}[theorem]{Remark}
\newtheorem{mlemma}[theorem]{Main Lemma}
\newtheorem{claim}[theorem]{Claim}

\newcommand{\alp}{\alpha}

\newcommand{\veps}{\varepsilon}
\newcommand{\del}{\delta}
\newcommand{\Del}{\Delta}
\newcommand{\ome}{\omega}
\newcommand{\Ome}{\Omega}
\newcommand{\bet}{\beta}
\newcommand{\gam}{\gamma}

\newcommand{\kap}{\kappa}
\newcommand{\sig}{\sigma}
\newcommand{\Sig}{\Sigma}

\newcommand{\Tht}{\Theta}
\newcommand{\lam}{\lambda}
\newcommand{\Lam}{\Lambda}
\newcommand{\vphi}{\varphi}

\newcommand{\fal}{\forall}
\newcommand{\exi}{\exists}

\newcommand{\Rarw }{\Rightarrow}

\newcommand{\lrarw}{\leftrightarrow}
\newcommand{\Lrarw}{\Leftrightarrow}

\newcommand{\cala}{{\cal A}}

\newcommand{\calc}{{\cal C}}

\newcommand{\calE}{{\cal E}}

\newcommand{\calg}{{\cal G}}
\newcommand{\calh}{{\cal H}}

\newcommand{\cals}{{\cal S}}

\newcommand{\calw}{{\cal W}}
\newcommand{\calx}{{\cal X}}
\newcommand{\caly}{{\cal Y}}

\newcommand{\la}{\langle}
\newcommand{\ra}{\rangle}

\title{
Well-foundedness proof for first-order reflection
}
\author{Toshiyasu Arai
\\
Graduate School of Science,
Chiba University
\\
1-33, Yayoi-cho, Inage-ku,
Chiba, 263-8522, JAPAN
\\
tosarai@faculty.chiba-u.jp
}
\date{}
\begin{document}
\maketitle

\section{Introduction}

In this note we show the following Theorem \ref{th:wf}.

\bth\label{th:wf}
${\sf KP}\Pi_{N}$ proves that {\rm each} initial segment 
$\{\alp\in OT: \alp<\psi_{\Ome}(\ome_{n}(\mK+1))\}$ {\rm is well-founded.}
\end{theorem}

{\rm KP}$\ell$ {\rm denotes a set theory for limits of admissibles.}
 {\rm KP}$\Pi_{N}$ {\rm denotes a set theory for $\Pi_{N}$-reflecting universes.}

$\Ome=\Ome_{1}=\ome_{1}^{CK}$, $\Ome_{\alp+1}=(\Ome_{\alp})^{+}$ for the next admissible above $\Ome_{\alp}$,
and $\Ome_{\lam}=\sup\{\Ome_{\alp}: 0<\alp<\lam\}$ for limit ordinals $\lam$.

$\alp=_{NF}\alp_{m}+\cdots+\alp_{0}$ means that $\alp=\alp_{m}+\cdots+\alp_{0}$ 
and $\alp_{m}\geq\cdots\geq\alp_{0}$
and each $\alp_{i}$ is a non-zero additive principal number.
For the binary Veblen function $\vphi$,
$\alp=_{NF}\vphi\bet\gam$ means that $\alp=\vphi\bet\gam$ and $\bet,\gam<\alp$.
$\alp=_{NF}\ome^{\bet}$ means that $\alp=\ome^{\bet}>\bet$.
$\alp=_{NF}\Ome_{\bet}$ means that $\alp=\Ome_{\bet}>\bet$.

Let $X<\alp :\Lrarw \fal\bet\in X(\bet<\alp)$, 
$\alp\leq X:\Lrarw \exi\bet\in X(\alp\leq\bet)$ and $X\leq Y:\Lrarw \fal \alp\in X\exi \bet\in Y(\alp\leq\bet)$.

IH denotes the Induction Hypothesis, MIH the Main IH and SIH the Subsidiary IH.

\section{Computable notation system $OT$}\label{subsec:decidable}

Let $\vec{\xi}=(\xi_{0},\ldots,\xi_{m-1})$ be a sequence of ordinals.
 The \textit{length} $lh(\vec{\xi}):=m$.
Sequences consisting of a single element $(\xi)$ is identified with the ordinal $\xi$,
and $\emptyset$ denotes the \textit{empty sequence}.
$\vec{0}$ denotes ambiguously a zero-sequence $(0,\ldots,0)$
with its length
$0\leq lh(\vec{0})\leq N-1$.
$\vec{\xi}*\vec{\mu}=(\xi_{0},\ldots,\xi_{m-1})*(\mu_{0},\ldots,\nu_{n-1})=(\xi_{0},\ldots,\xi_{m-1},\mu_{0},\ldots,\mu_{n-1})$
denotes the \textit{concatenated} sequence of $\vec{\xi}$ and $\vec{\mu}$.


$
\Lam=\veps(\mK)=\veps_{\mK+1}
$
denotes the next epsilon number above the least $\Pi_{N-2}$-
\\
indescribable cardinal
$\mK$, and
$\veps(\Lam)=\veps_{\mK+2}$ the next epsilon number above $\Lam$.

\bdf\label{df:Lam}
{\rm

For a non-zero ordinal $\xi<\veps(\Lam)$, its Cantor normal form with base $\Lam$
is uniquely determined as
\beqn\label{eq:CantornfLam}
\xi=_{NF}\sum_{i\leq m}\Lam^{\xi_{i}}a_{i}=\Lam^{\xi_{m}}a_{m}+\cdots+\Lam^{\xi_{0}}a_{0}
\eeqn
where
$\xi_{m}>\cdots>\xi_{0}, \, 0<a_{i}<\Lam$.

  \benu
  \item\label{df:Lam2.0}
  $K(\xi)=\{a_{i}:i\leq m\}\cup\bigcup\{K(\xi_{i}):i\leq m\}$ is the set of \textit{components} of $\xi$
  with $K(0)=\emptyset$.
For a sequence $\vec{\xi}=(\xi_{0},\ldots,\xi_{n-1})$ of ordinals $\xi_{i}<\veps(\Lam)$,
$K(\vec{\xi})
:=\bigcup\{K(\xi_{i}): i<n\}$.

 \item\label{df:Exp2.1}
 For $\xi>1$,
 $te(\xi)=\xi_{0}$ in (\ref{eq:CantornfLam}) is the \textit{tail exponent}, 
 and
 $he(\xi)=\xi_{m}$ is the \textit{head exponent} of $\xi$, resp.
The \textit{head} $Hd(\xi):=\Lam^{\xi_{m}}a_{m}$, and the \textit{tail}
 $Tl(\xi):=\Lam^{\xi_{0}}a_{0}$ of $\xi$.


 \item\label{df:Exp2.2h}
 $he^{(i)}(\xi)$ is the \textit{$i$-th head exponent} of $\xi$, defined recursively by
 \\
 $he^{(0)}(\xi)=\xi$, $he^{(i+1)}(\xi)=he(he^{(i)}(\xi))$.
 
 The \textit{$i$-th tail exponent} $te^{(i)}(\xi)$ is defined similarly.
 
 \item\label{df:Exp2.3}
 $\zeta$ is a \textit{part} of $\xi$, denoted by
 $\zeta\leq_{pt}\xi$ iff
 \\
 $\zeta=_{NF}\sum_{i\geq n}\Lam^{\xi_{i}}a_{i}=\Lam^{\xi_{m}}a_{m}+\cdots+\Lam^{\xi_{n}}a_{n}$ {\rm for an} $n\, (0\leq n\leq m+1)$.
 
$\zeta<_{pt}\xi:\Lrarw\zeta\leq_{pt}\xi \spand \zeta\neq\xi$.
 
\item\label{df:Exp2.4}
A sequence $\vec{\mu}=(\mu_{0},\ldots,\mu_{n})$ is an \textit{iterated tail parts} of $\xi$, denoted by
$\vec{\mu}\subset_{pt}\xi$
iff $ \mu_{0}\leq_{pt}\xi \spand \fal i<n(\mu_{i+1}\leq_{pt} te(\mu_{i}))$.

\item\label{df:Exp2.5}
$\vec{\nu}=(\nu_{0},\ldots,\nu_{n})*\vec{0}<\xi$ iff there exists a sequence
$\vec{\mu}=(\mu_{0},\ldots,\mu_{n})$ such that
$\vec{\mu}\subset_{pt}\xi$ and $\nu_{i}<\mu_{i}$ for every $i\leq n$.

\item\label{df:Exp2.8}
Let $\vec{\nu}=(\nu_{0},\ldots,\nu_{n})$ and $\vec{\xi}=(\xi_{0},\ldots,\xi_{n})$ be
sequences of ordinals in the same length, and $0\leq k\leq n$.

$\vec{\nu}<_{k}\vec{\xi}:\Lrarw \fal i<k(\nu_{i}\leq\xi_{i})\land 
(\nu_{k},\ldots,\nu_{n})< \xi_{k}$.

\item\label{df:Exp2.9}
 $\zeta$ is a \textit{step-down} of $\xi$, denoted by
 $\zeta<_{sd}\xi$ iff
 \\
 $\zeta=\Lam^{\xi_{m}}a_{m}+\cdots+\Lam^{\xi_{1}}a_{1}+\Lam^{\xi_{0}}b+\nu$ for some ordinals
 $b<a_{0}$ and $\nu<\Lam^{\xi_{0}}$.

\item\label{df:Exp2.10}
$\vec{\nu}=(\nu_{0},\ldots,\nu_{n})*\vec{0}<_{sd}\xi$ iff 
$\nu_{i}<_{sd}te^{(i)}(\xi)$ for every $i\leq n$.

\item\label{df:Exp2.11}
$\zeta\leq_{sp}\xi:\Lrarw\exi\mu\leq_{pt}\xi(\zeta\leq_{sd}\mu)$, and
$\zeta<_{sp}\xi:\Lrarw\exi\mu\leq_{pt}\xi(\zeta<_{sd}\mu)$.

\item\label{df:Exp2.12}
$\vec{\nu}<_{sp}\xi$ iff 
$\vec{\nu}<_{sd}\mu$ for a $\mu\leq_{pt}\xi$.

Let $p(\vec{\nu},\xi)$ denote the number $p\, (0\leq p<m)$ such that
$\xi=_{NF}\mu+\sum_{i<p}\Lambda^{\xi_{i}}a_{i}$
for
$\mu=\Lambda^{\xi_{m}}a_{m}+\cdots+\Lambda^{\xi_{p}}a_{p}$ and $\vec{\nu}<_{sd}\mu$.
\eenu

}
\edf

It is easy to see that $\vec{\nu}<_{sd}\xi \Rarw \vec{\nu}<\xi$.

\bprp\label{prp:idless}
$\vec{\nu}<\xi\leq\zeta \Rarw \vec{\nu}<\zeta$.
\eprp

\bdf\label{df:nfform}
{\rm A sequence of ordinals} $\vec{\xi}=(\xi_{2},\ldots,\xi_{N-1})$ {\rm is said to be} irreducible 
{\rm iff}
$\fal i<N-1\fal k>0(\xi_{i}>0 \Rarw Tl(\xi_{i})\geq\Lam_{k}(\xi_{i+k}+1))$.
 \edf
 
 \bdf\label{df:lx}
 {\rm Let} $\vec{\nu}\neq\vec{\xi}$ {\rm and let} $i$ {\rm be the minimal number such that}
 $\nu_{i}\neq\xi_{i}$. 
 {\rm Suppose} $(\xi_{i},\ldots,\xi_{N-1})\neq\vec{0}${\rm , and let}
  $k_{1}\geq i$ {\rm be the minimal number such that} $\xi_{k_{1}}\neq 0$. {\rm Then}
  $\vec{\nu}<_{lx,2}\vec{\xi}$ {\rm iff one of the followings holds:}
  \benu
  \item
  $(\nu_{i},\ldots,\nu_{N-1})=\vec{0}$.
  
  \item
  {\rm In what follows assume}  $(\nu_{i},\ldots,\nu_{N-1})\neq\vec{0}${\rm , and let}
  $k_{0}\geq i$ {\rm be the minimal number such that} $\nu_{k_{0}}\neq 0\, (i=\min\{k_{0},k_{1}\})$. {\rm Then}
  $\vec{\nu}<_{lx,2}\vec{\xi}$ {\rm iff one of the followings holds:}
  
   \benu
   \item
   $i=k_{0}< k_{1}$ {\rm and} $he^{(k_{1}-i)}(\nu_{i})\leq\xi_{k_{1}}$.
   \item
   $k_{0}\geq k_{1}=i$ {\rm and} $\nu_{k_{0}}<he^{(k_{0}-i)}(\xi_{i})$.
   \eenu
  \eenu
\edf

\bdf\label{df:SD}
{\rm
A set $SD$ of sequences $\vec{\xi}=(\xi_{2},\ldots,\xi_{N-1})$ of ordinals $\xi_{i}<\veps(\Lam)$
is defined recursively as follows.
\benu
\item
$\vec{0}*(a)\in SD$ for each $a<\Lam$.
\item
(Cf.\,Definition \ref{df:Lam}.\ref{df:Exp2.10}.)
Let $\vec{\xi}=(\xi_{2},\ldots,\xi_{N-1})\in SD$, $1\leq k<N-1$, $\zeta<\veps(\Lam)$ be an ordinal such that 
$(\xi_{k+1},\ldots,\xi_{N-1})<_{sd}\zeta$, and
$(\xi_{2},\ldots,\xi_{k-1},\xi_{k},\zeta)*\vec{0}\in SD$.
Then for
 $\zeta_{k}=\xi_{k}+\Lam^{\zeta}a$ with an ordinal $a<\Lam$,
$(\xi_{2},\ldots,\xi_{k-1})*(\zeta_{k})*(\xi_{k+1},\ldots,\xi_{N-1})\in SD$
 and $(\xi_{2},\ldots,\xi_{k-1})*(\zeta_{k})*\vec{0}\in SD$.
\eenu
}
\edf

\bprp\label{prp:SD}
Let $\vec{\xi}=(\xi_{2},\ldots,\xi_{N-1})\in SD$.
\benu
\item\label{prp:SD.-1}
$(\xi_{2},\ldots,\xi_{i})*\vec{0}\in SD$ for each $i$ with $1\leq i<N$.
\item\label{prp:SD.0}
For $2\leq i<j<k<N$, if $\xi_{i}\neq 0$ and $\xi_{k}\neq 0$, then $\xi_{j}\neq 0$.

\item\label{prp:SD.1}
Let $\xi_{i}\neq 0$.
Then 
$(\xi_{i+1},\ldots,\xi_{N-1})<_{sd}te(\xi_{i})$.
\item\label{prp:SD.2}
$\vec{\xi}$ is irreducible.
\eenu
\eprp

\bdf\label{df:lxo}
{\rm Let us define an ordinal $o(\vec{\xi})$ for irreducible $\vec{\xi}=(\xi_{2},\ldots,\xi_{N-1})$ by
\[
o(\vec{\xi})=\sum\{\Lam_{i-1}(\xi_{i}+1): 2\leq i\leq N-1, \xi_{i}\neq 0\}
\]
In particular $o(\vec{0})=0$.}
\edf
Note that we have $Tl(\xi_{i})\geq\Lam_{k}(\xi_{i+k}+1)$ for $\xi_{i}\neq 0$ and irreducible $\vec{\xi}$.
Therefore $\xi_{i}+\Lam_{k}(\xi_{i+k}+1)=\xi_{i}\#\Lam_{k}(\xi_{i+k}+1)$ for the natural sum $\#$.

\bprp\label{prp:lxo}
For irreducible $\vec{\nu},\vec{\xi}$,
\[
\vec{\nu}<_{lx,2}\vec{\xi} \Rarw o(\vec{\nu})<o(\vec{\xi})
.\]
\eprp
\bprf
Let $\vec{\nu}<_{lx,2}\vec{\xi}$. 
Then $\vec{\nu}\neq\vec{\xi}$ and let $i$ be the minimal number such that
 $\nu_{i}\neq\xi_{i}$. 
It suffices to show that 
$a_{0}=o((\nu_{i},\ldots,\nu_{N-1}))<o((\xi_{i},\ldots,\xi_{N-1}))=a_{1}$, where 
$o((\xi_{i},\ldots,\xi_{N-1}))=\sum\{\Lam_{j-1}(\xi_{j}+1): i\leq j\leq N-1, \xi_{j}\neq 0\}$.

 We have $(\xi_{i},\ldots,\xi_{N-1})\neq\vec{0}$, and let
 $k_{1}\geq i$ be the minimal number such that $\xi_{k_{1}}\neq 0$. 
 When $(\nu_{i},\ldots,\nu_{N-1})\neq\vec{0}$, let
  $k_{0}\geq i$ be the minimal number such that $\nu_{k_{0}}\neq 0$. 
 One of the following cases occurs, cf. Definition \ref{df:lx}.
 \\
{\bf Case 0}.
  $(\nu_{i},\ldots,\nu_{N-1})=\vec{0}$: Then $a_{0}=0<a_{1}$.
\\
{\bf Case 1}.
   $i=k_{0}< k_{1}=i+k$ {\rm and} $he^{(k)}(\nu_{i})\leq\xi_{i+k}$:
We have by $k>0$
\\
$o((\nu_{i},\ldots,\nu_{N-1}))=\Lam_{i-1}(\nu_{i}+1)+o((\nu_{i+1},\ldots,\nu_{N-1}))<\Lam_{i+k-1}(he^{(k)}(\nu_{i})+1)$.
On the other hand we have 
$o((\xi_{i},\ldots,\xi_{N-1}))=o((0,\ldots,0,\xi_{i+k},\ldots,\xi_{N-1}))=\Lam_{i+k-1}(\xi_{i+k}+1)+o(\xi_{i+k+1},\ldots,\xi_{N-1})\geq\Lam_{i+k-1}(\xi_{i+k}+1)$.
Hence $a_{0}<a_{1}$.
\\
{\bf Case 2}.
   $i+k=k_{0}\geq k_{1}=i$ {\rm and} $\nu_{i+k}<he^{(k)}(\xi_{i})$: Then 
\beqnarrs
o((\nu_{i},\ldots,\nu_{N-1})) & = & o((0,\ldots,0,\nu_{i+k},\ldots,\nu_{N-1}))
\\
 & = &  \Lam_{i+k-1}(\nu_{i+k}+1)+o((\nu_{i+k+1},\ldots,\nu_{N-1}))
 \\
 & < & \Lam_{i+k-1}(\nu_{i+k}+1)\cdot 2\leq \Lam_{i+k-1}(he^{(k)}(\xi_{i}))\cdot 2
\eeqnarrs
On the other hand we have by $i>1$ and $\xi_{i}\geq\Lam_{k}(he^{(k)}(\xi_{i}))$
\beqnarrs
o((\xi_{i},\ldots,\xi_{N-1})) & = & \Lam_{i-1}(\xi_{i}+1)+o((\xi_{i+1},\ldots,\xi_{N-1}))
\\
& \geq & \Lam_{i-1}(\xi_{i}+1) >
\Lam_{i+k-1}(he^{(k)}(\xi_{i}))\cdot 2
\eeqnarrs
Hence $a_{0}<a_{1}$.
\eprf
\\

{\rm The sets of ordinal terms} $OT\subset\Lam=\veps_{\mK+1}$ {\rm and} $E\subset\veps_{\mK+2}$
{\rm over symbols} $\{0,\mK,\Lam,+,\ome,\vphi, \Ome, \psi\}$
{\rm together with}
 {\rm sequences} 
$\vec{m}(\alp)=(m_{k}(\alp))_{2\leq k\leq N-1}\, (\alp\in OT\cap\mK)$, 
 and
finite sets $K_{\del}(\alp)\subset OT$ for $\alp\in OT$
{\rm are defined by simultaneous recursion as follows.}

Let $pd(\psi_{\pi}^{\vec{\nu}}(a))=\pi$ (even if $\vec{\nu}=\vec{0}$).
Moreover for $n$,
$pd^{(n)}(\alp)$ is defined recursively by $pd^{(0)}(\alp)=\alp$ and
$pd^{(n+1)}(\alp)\simeq pd(pd^{(n)}(\alp))$.

For terms $\pi,\kap\in OT$,
$\pi\prec\kap$ denotes the transitive closure of the relation
$\{(\pi,\kap): \exi \vec{\xi}\exi b[\pi=\psi_{\kap}^{\vec{\xi}}(b)]\}$,
 and its reflexive closure
$\pi\preceq\kap:\Lrarw \pi\prec\kap \lor \pi=\kap\Lrarw\exi n(\kap=pd^{(n)}(\pi))$.

For each ordinal term $\alp=\psi_{\pi}^{\vec{\nu}}(a)$,
a series $(\pi_{i})_{i\leq L}$ of ordinal terms is uniquely determined as follows:
$\pi_{L}=\alp$, $\pi_{i}=pd(\pi_{i+1})$ and $\pi_{0}=\mK$.
Let us call the series $(\pi_{i})_{i\leq L}$ the \textit{collapsing series} of $\alp=\pi_{L}$.

Then we see that an ordinal term
$\alp=\psi_{\pi}^{\vec{\nu}}(a)$ with $\vec{\nu}\neq\vec{0}$ 
is constructed by Definition \ref{df:notationsystem}.\ref{df:notationsystem.10} below
iff $L=1$.
$\alp$ is constructed by Definition \ref{df:notationsystem}.\ref{df:notationsystem.12}
iff $L\equiv 1 \!\!\!\!\pmod{(N-2)}$.
Otherwise $\alp$ is constructed by Definition \ref{df:notationsystem}.\ref{df:notationsystem.11}.

\bdf\label{df:notationsystem}
{\rm
$\ell\alp$ denotes the number of occurrences of symbols
\\
$\{0,\mK,\Lam,+,\ome,\vphi, \Ome, \psi\}$
in terms $\alp\in OT\cup E$.

\benu
\item\label{df:notationsystem.1}

 \benu
 \item\label{df:notationsystem.3}
$0\in E$.

 \item\label{df:notationsystem.5-1}
If $0<a\in OT$,
then $a\in E$.
$K(a)=\{a\}$.

 \item\label{df:notationsystem.5}
If $\{\xi_{i}:i\leq m\}\subset E$, $\xi_{m}>\cdots>\xi_{0}>0$ and
$0<b_{i}\in OT$,
then $\sum_{i\leq m}\Lam^{\xi_{i}}b_{i}=
\Lam^{\xi_{m}}b_{m}+\cdots+\Lam^{\xi_{0}}b_{0}\in E$.
$K(\sum_{i\leq m}\Lam^{\xi_{i}}b_{i})=\{b_{i}:i\leq m\}\cup\bigcup\{K(\xi_{i}):i\leq m\}$.

\item
For sequences $\vec{\nu}=(\nu_{2},\ldots,\nu_{N-1})$, let
$K(\vec{\nu})=\bigcup_{2\leq i\leq N-1}K(\nu_{i})$.

  
 

 \eenu
 
 \item
 \benu
 \item\label{df:notationsystem.2}
$0,\mK\in OT$.
$m_{k}(0)=0$ for any $k$, and $K_{\del}(0)=K_{\del}(\mK)=\emptyset$.

 \item\label{df:notationsystem.4}
If $\alp=_{NF}\alp_{m}+\cdots+\alp_{0}\, (m>0)$ with $\{\alp_{i}:i\leq m\}\subset OT$, 
then
$\alp\in OT$, and $m_{k}(\alp)=0$ for any $k$.
$K_{\del}(\alp)=K_{\del}(\alp_{0},\ldots,\alp_{m})$.

 \item\label{df:notationsystem.6}
If $\alp=_{NF}\vphi\bet\gam$ with $\{\bet,\gam\}\subset OT\cap\mK$, then
$\alp\in OT$, and $m_{k}(\alp)=0$ for any $k$.
$K_{\del}(\alp)=K_{\del}(\bet,\gam)$.

 \item\label{df:notationsystem.7}
If $\alp=_{NF}\ome^{\bet}$ with $\mK<\bet\in OT$, then $\alp\in OT$,
and $m_{k}(\alp)=0$ for any $k$.
$K_{\del}(\alp)=K_{\del}(\bet)$.

 \item\label{df:notationsystem.8}
If $\alp=_{NF}\Ome_{\bet}$ with $\bet\in OT\cap\mK$, then
$\alp\in OT$.
$m_{2}(\alp)=1, m_{k}(\alp)=0$ for any $k>2$ if $\bet$ is a successor ordinal.
Otherwise
$m_{k}(\alp)=0$ for any $k$.
In each case
$K_{\del}(\alp)=K_{\del}(\bet)$.

 \item\label{df:notationsystem.9}
Let $\alp=\psi_{\pi}(a):=\psi_{\pi}^{\vec{0}}(a)$ where $\pi$ is a regular term , i.e.,
either $\pi=\mK$ or
$\vec{m}(\pi)\neq\vec{0}$,
and  $K_{\alp}(\pi,a)<a$.

Then
$\alp=\psi_{\pi}(a)\in OT$. Let
$m_{k}(\alp)=0$ for any $k$.
$K_{\del}(\psi_{\pi}(a))=\emptyset$ if $\alp<\del$.
$K_{\del}(\psi_{\pi}(a))=\{a\}\cup K_{\del}(a,\pi)$ otherwise.

 \item\label{df:notationsystem.10}
Let $\alp=\psi_{\mK}^{\vec{\nu}}(a)$ with
 $\vec{\nu}=\vec{0}*(b)\, (lh(\vec{\nu})=N-2)$ and $b,a\in OT$
such that $0<b\leq a$ and $K_{\alp}(b,a)<a$.

 Then
$\alp=\psi_{\mK}^{\vec{\nu}}(a)\in OT$.
Let
$m_{N-1}(\alp)=b$, $m_{k}(\alp)=0$
for $k<N-1$.
$K_{\del}(\psi_{\mK}^{\vec{\nu}}(a))=\emptyset$ if $\alp<\del$.
$K_{\del}(\psi_{\mK}^{\vec{\nu}}(a))=\{a\}\cup \bigcup\{K_{\del}(\gam):\gam\in K(\nu)\}$ otherwise.

 \item\label{df:notationsystem.11}
Let $\pi\in OT\cap\mK$ be such that 
$m_{k+1}(\pi)\neq 0$ 
and $\fal i>k+1(m_{i}(\pi)=0)$
for a $k\, (2\leq k\leq N-2)$, and
$b,a\in OT$ such that $0< b\leq a$.
Let $\vec{\nu}=(\nu_{2},\ldots,\nu_{N-1})$ be a sequence
defined by $\fal i<k(\nu_{i}=m_{i}(\pi))$,
$\nu_{k}=m_{k}(\pi)+\Lam^{m_{k+1}(\pi)}b$,
and $\fal i>k(\nu_{i}=0)$.

Then 
$\alp=\psi_{\pi}^{\vec{\nu}}(a)\in OT$ if
$K_{\alp}(\pi,a,b)\cup K_{\alp}(K(\vec{m}(\pi)))<a$.
Let $m_{i}(\alp)=\nu_{i}$ for each $i$.
$K_{\del}(\psi_{\pi}^{\vec{\nu}}(a))=\emptyset$ if $\alp<\del$.
Otherwise
$K_{\del}(\psi_{\pi}^{\vec{\nu}}(a))=
\{a\}\cup K_{\del}(a,\pi)\cup\bigcup\{K_{\del}(b): b\in K(\vec{\nu})\}$.

 \item\label{df:notationsystem.12}
Let $\pi\in OT\cap\mK$ be such that 
$m_{2}(\pi)\neq 0$ and $\fal i>2(m_{i}(\pi)=0)$,
and $a\in OT$.
Let $\vec{0}\neq\vec{\nu}=(\nu_{2},\ldots,\nu_{N-1})\in SD$ be a sequence
 of ordinal terms $\nu_{i}\in E$
such that $\vec{\nu}<_{sp} m_{2}(\pi)$.


Then
$\alp=\psi_{\pi}^{\vec{\nu}}(a)$ if $K_{\alp}(\pi,a)<a$, and
\beqn\label{eq:notationsystem.12}
\forall k( K_{\alpha}(\nu_{k})<\max K(\nu_{k}))
\eeqn

Let $m_{i}(\alp)=\nu_{i}$ for each $i$.

$K_{\del}(\psi_{\pi}^{\vec{\nu}}(a))=\emptyset$ if $\alp<\del$.
Otherwise
$K_{\del}(\psi_{\pi}^{\vec{\nu}}(a))=\{a\}\cup K_{\del}(a,\pi)\cup\bigcup\{K_{\del}(b): b\in K(\vec{\nu})\}$.

 \eenu

\eenu
}
\edf



Let $\calh_{\gam}(\del)$ be the Skolem hull defined in \cite{KPpiNsmpl}
such that
for any $\alp\in OT$ and any $\del$ such that $\del=0,\Lam$ or $\del=\psi_{\pi}^{\vec{\nu}}(b)$ for some $\pi,b,\vec{\nu}$,
$\alp\in\calh_{\gam}(\del) \Lrarw K_{\del}(\alp)<\gam$.

\bprp\label{prp:psicomparison}{\rm (Cf. Proposition 2.19 in \cite{KPpiNsmpl})}\\
Let $\vec{\nu}=(\nu_{0},\ldots,\nu_{N-3})$, $\vec{\xi}=(\xi_{0},\ldots,\xi_{N-3})$ be 
irreducible sequences of ordinals$<\veps_{\Lam+2}$, and assume that
$\psi_{\pi}^{\vec{\nu}}(b)<\pi$ and $\psi_{\kap}^{\vec{\xi}}(a)<\kap$.

Then $\bet_{1}=\psi_{\pi}^{\vec{\nu}}(b)<\psi_{\kap}^{\vec{\xi}}(a)=\alp_{1}$ iff one of the following cases holds:
\benu
\item\label{prp:psicomparison.0}
$\pi\leq \psi_{\kap}^{\vec{\xi}}(a)$.

\item\label{prp:psicomparison.1}
$b<a$, $\psi_{\pi}^{\vec{\nu}}(b)<\kap$ and 
$K(\vec{\nu})\cup\{\pi,b\}\subset\calh_{a}(\psi_{\kap}^{\vec{\xi}}(a))$.

\item\label{prp:psicomparison.2}
$b>a$ and $K(\vec{\xi})\cup\{\kap,a\}\not\subset\calh_{b}(\psi_{\pi}^{\vec{\nu}}(b))$.

\item\label{prp:psicomparison.25}
$b=a$, $\kap<\pi$ and $\kap\not\in\calh_{b}(\psi_{\pi}^{\vec{\nu}}(b))$.

\item\label{prp:psicomparison.3}
$b=a$, $\pi=\kap$, $K(\vec{\nu})\subset\calh_{a}(\psi_{\kap}^{\vec{\xi}}(a))$, and
$\vec{\nu}<_{lx,2}\vec{\xi}$.

\item\label{prp:psicomparison.4}
$b=a$, $\pi=\kap$, 
$K(\vec{\xi})\not\subset\calh_{b}(\psi_{\pi}^{\vec{\nu}}(b))$.

\eenu

\eprp

\bprp\label{prp:G6}
\benu
\item\label{prp:G6.1}
$\alp\leq\bet \Rarw K_{\alp}(\gam)\supset K_{\bet}(\gam)$.


\item\label{prp:G6.3}
Let $\bet=\psi_{\pi}^{\vec{\nu}}(b)$ with $\pi=\psi_{\kap}^{\vec{\xi}}(a)$.
Then $a<b$.

\item\label{prp:G6.4}
If $\kap<\psi_{\pi}^{\vec{\nu}}(b)<\kap^{+}$, then $\pi=\kap^{+}$(, and $\vec{\nu}=\emptyset$).

\item\label{prp:G6.5}
For $\alpha=\psi_{\pi}^{\vec{\nu}}(a)\in OT$,
$\max K(\vec{\nu})\leq a$ holds.
\eenu
\eprp
\bprf
\ref{prp:G6}.\ref{prp:G6.1} is seen by induction on $\ell\gam$.
\\

\noindent
\ref{prp:G6}.\ref{prp:G6.3}.
Let $\bet=\psi_{\pi}^{\vec{\nu}}(b)$ with $\pi=\psi_{\kap}^{\vec{\xi}}(a)$.
Then $K_{\bet}(\{\pi,b\}\cup K(\vec{\nu}))<b$.
On the other hand we have $\bet<\pi$.
Hence $a\in K_{\bet}(\pi)<b$.
\\

\noindent
\ref{prp:G6}.\ref{prp:G6.4}.
Let $\kap<\psi_{\pi}^{\vec{\nu}}(b)<\kap^{+}$.
If $\vec{\nu}\neq\emptyset$, then $\kap^{+}<\psi_{\pi}^{\vec{\nu}}(b)$.
Hence $\vec{\nu}=\emptyset$.
Let $\kap=\Ome_{a}\geq a$ with $\kap^{+}=\Ome_{a+1}$. 
Then $a\in\calh_{b}(\psi_{\pi}(b))$, and $\Ome_{a+1}\in\calh_{b}(\psi_{\pi}(b))$.
If $\kap^{+}=\Ome_{a+1}<\pi$, then $\kap^{+}<\psi_{\pi}(b)$.
Hence $\kap<\pi\leq\kap^{+}$, and $\pi=\kap^{+}$.
\\

\noindent
\ref{prp:G6}.\ref{prp:G6.5}.
This is seen by induction on $\ell\alpha$.
Ww have $c<a$ by Proposition \ref{prp:G6}.\ref{prp:G6.3} when $\pi=\psi_{\sig}^{\vec{\mu}}(c)$

When $\alpha$ is constructed by Definition \ref{df:notationsystem}.\ref{df:notationsystem.11},
$\nu_{k}=m_{k}(\pi)+\Lam^{m_{k+1}(\pi)}b$ holds for $b\leq a$.
By IH we have $\max K(\vec{m}(\pi))\leq c<a$ when $\pi=\psi_{\sig}^{\vec{\mu}}(c)$.

Suppose $\alpha$ is constructed by Definition \ref{df:notationsystem}.\ref{df:notationsystem.12}.
We obtain $\vec{\nu}<_{sp}m_{2}(\pi)$, and hence
$\max K(\vec{\nu})\leq \max K(m_{2}(\pi))\leq c<a$ by IH.
\eprf
\\



Let $OT_{n}$ denote the subsystem of $OT$ such that $\alp\in OT_{n}$ iff each ordinal subterm occurring in $\alp$ is smaller than
$\ome_{n}(\mK+1)$.

\bdf\label{df:notationsystemnslice}

\benu

\item\label{df:notationsystemnslice.2}
{\rm For} $\alp\in OT\cap \mK$, $\alp\in OT_{n} \Lrarw \calE(\alp)\subset OT_{n}$.

\item\label{df:notationsystemnslice.7}
{\rm If} $\alp=_{NF}\ome^{\bet}<\ome_{n}(\mK+1)$ {\rm with} $\mK<\bet\in OT_{n}$, {\rm then} $\alp\in OT_{n}$,

\item\label{df:notationsystemnslice.11}
{\rm Let} $\alp=\psi_{\pi}^{\vec{\nu}}(a)\in OT$ {\rm such that} $\{a,\pi\}\cup K(\vec{\nu})\subset OT_{n}$.
 {\rm Then} 
$\alp=\psi_{\pi}^{\vec{\nu}}(a)\in OT_{n}$.

\eenu

\edf


\bprp\label{lem:Tn}

For any $n<\ome$ and $\del=\psi_{\Ome}(\ome_{n}(\mK+1))$,
\benu
\item\label{lem:Tn1}
$\fal\alp\in OT\cap\psi_{\Ome}(\ome_{n}(\mK+1)) (\alp\in OT_{n})$.
\item\label{lem:Tn2}
$\fal\alp\in OT(\{\alp\}\cup K_{\del}(\alp)<\ome_{n}(\mK+1) \Rarw \alp\in OT_{n})$.
\eenu
\eprp
\bprf
These are shown simultaneously by induction on $\ell\alp$ for $\alp\in OT$.

If $\alp$ is not a strongly critical number, then IH yields the lemmas.
Let $\alp=\psi_{\pi}^{\vec{\nu}}(b)$ for some $\pi,\vec{\nu},b$.
\\

\noindent
\ref{lem:Tn}.\ref{lem:Tn1}.
Let $\alp<\psi_{\Ome}(\ome_{n}(\mK+1))$.
Since $\Ome$ is the least recursively regular ordinal, $\psi_{\Ome}(\ome_{n}(\mK+1))<\Ome$
 and $K_{\alp}(\{\Ome,\ome_{n}(\mK+1)\})=\emptyset$,
we see that 
$b<\ome_{n}(\mK+1)$, $\psi_{\pi}^{\vec{\nu}}(b)<\Ome$ and 
$K_{\del}(K(\vec{\nu})\cup\{\pi,b\})<\ome_{n}(\mK+1)$.
By $\psi_{\pi}^{\vec{\nu}}(b)<\Ome$ we obtain $\pi=\Ome$, and $\vec{\nu}=\emptyset$.
IH on Proposition \ref{lem:Tn}.\ref{lem:Tn1} with $\{b\}\cup K_{\del}(b)<\ome_{n}(\mK+1)$ yields $b\in OT_{n}$.
\\

\noindent
\ref{lem:Tn}.\ref{lem:Tn2}.
Let $K_{\del}(\alp)<\ome_{n}(\mK+1)$.
If $\alp<\del=\psi_{\Ome}(\ome_{n}(\mK+1))$, then $\alp\in OT_{n}$ by Proposition \ref{lem:Tn}.\ref{lem:Tn1}.
Suppose $\alp\geq\del$.
Then $K_{\del}(\alp)=\{b\}\cup K_{\del}(K(\vec{\nu})\cup\{\pi,b\})$.
IH with $\pi\leq\mK$ yields $\{b,\pi\}\subset OT_{n}$.
In particular $b<\ome_{n}(\mK+1)$.
This yields $K(\vec{\nu})\leq b<\ome_{n}(\mK+1)$ by
Definition \ref{df:notationsystem}.\ref{df:notationsystem.10} and \ref{df:notationsystem}.\ref{df:notationsystem.11}.
Hence by IH we obtain $K(\vec{\nu})\subset OT_{n}$.
\eprf
\\

\noindent 
Therefore it suffices show the following Theorem \ref{lem:wfTn} to prove Theorem \ref{th:wf}.

\bth\label{lem:wfTn} 
For {\rm each} $n<\ome$,
${\sf KP}\Pi_{N}$ proves that $(OT_{n},<)$ is well-founded.
\end{theorem}

\subsection{Coefficients}
In this subsection we introduce coefficient sets $\calE(\alp),G_{\kap}(\alp), F_{\del}(\alp),k_{\del}(\alp)$ of $\alp\in OT$,
each of which is a finite set of subterms of $\alp$.
These are utilized in our wellfoundedness proof in section \ref{sec:distinguished}.
Roughly $\calE(\alp)$ is the set of subterms of the form $\psi_{\pi}^{\vec{\nu}}(a)$, and
$F_{\del}(\alp)$ [$k_{\del}(\alp)$] the set of subterms$<\del$ [subterms$\geq\del$], resp.
$G_{\kap}(\alp)$ is an analogue of sets $K_{\kap}\alp$ in \cite{odMahlo}.

Let $pd(\psi_{\pi}^{\vec{\nu}}(a))=\pi$ (even if $\vec{\nu}=\emptyset$).
Moreover for $n$,
$pd^{(n)}(\alp)$ is defined recursively by $pd^{(0)}(\alp)=\alp$ and
$pd^{(n+1)}(\alp)\simeq pd(pd^{(n)}(\alp))$.

{\rm For terms} $\pi,\kap\in OT$,
$\pi\prec\kap$ {\rm denotes the transitive closure of the relation}
$\{(\pi,\kap): \exi \vec{\xi}\exi b[\pi=\psi_{\kap}^{\vec{\xi}}(b)]\}$,
{\rm and its reflexive closure}
$\pi\preceq\kap:\Lrarw \pi\prec\kap \lor \pi=\kap$.

For terms $\pi,\kap\in OT$,
$\pi\prec\kap$ denotes the transitive closure of the relation
$\{(\pi,\kap): \exi \vec{\xi}\exi b[\pi=\psi_{\kap}^{\vec{\xi}}(b)]\}$,
 and its reflexive closure
$\pi\preceq\kap:\Lrarw \pi\prec\kap \lor \pi=\kap$.

\bdf

{\rm For terms} $\alp,\kap,\del\in OT$, {\rm finite sets} $\calE(\alp),  G_{\kap}(\alp), F_{\del}(\alp), k_{\del}(\alp)\subset OT$ {\rm are defined recursively as follows.}
\benu

\item
$\calE(\alp)=\emptyset$ {\rm for} $\alp\in\{0,\mK\}$.
$\calE(\alp_{m}+\cdots+\alp_{0})=\bigcup_{i\leq m}\calE(\alp_{i})$.
$\calE(\vphi\bet\gam)=\calE(\bet)\cup\calE(\gam)$.
$\calE(\ome^{\alp})=\calE(\alp)$.
$\calE(\Ome_{\alp})=\calE(\alp)$.

\item
$\calE(\psi_{\pi}^{\vec{\nu}}(a))=\{\psi_{\pi}^{\vec{\nu}}(a)\}$.

\item
$\cala(\alp)=\bigcup\{\cala(\bet): \bet\in\calE(\alp)\}$
{\rm for} $\cala\in\{G_{\kap},F_{\del},k_{\del}\}$.

\item
\[
G_{\kap}(\psi_{\pi}^{\vec{\nu}}(a))=\left\{
\begin{array}{ll}
G_{\kap}(\{\pi,a\}\cup K(\vec{\nu})) & \kap<\pi
\\
G_{\kap}(\pi) & \pi<\kap \spand \pi\not\preceq\kap
\\
\{\psi_{\pi}^{\vec{\nu}}(a)\} & \pi\preceq\kap
\end{array}
\right.
\]
\[
F_{\del}(\psi_{\pi}^{\vec{\nu}}(a))=\left\{
\begin{array}{ll}
F_{\del}(\{\pi,a\}\cup K(\vec{\nu})) & \psi_{\pi}^{\vec{\nu}}(a)\geq\del
\\
\{\psi_{\pi}^{\vec{\nu}}(a)\} & \psi_{\pi}^{\vec{\nu}}(a)<\del
\end{array}
\right.
\]
\[
k_{\del}(\psi_{\pi}^{\vec{\nu}}(a))=\left\{
\begin{array}{ll}
\{\psi_{\pi}^{\vec{\nu}}(a)\}\cup k_{\del}(\{\pi,a\}\cup K(\vec{\nu})) & \psi_{\pi}^{\vec{\nu}}(a)\geq\del
\\
\emptyset & \psi_{\pi}^{\vec{\nu}}(a)<\del
\end{array}
\right.
\]
\eenu

{\rm For} $\cala\in\{K_{\del},G_{\kap},F_{\del},k_{\del}\}$ {\rm and sets} $X\subset OT$,
$\cala(X):=\bigcup\{\cala(\alp): \alp\in X\}$.
\edf

\bdf
{\rm $S(\eta)$ denotes the set of immediate subterms of $\eta$ when $\eta\not\in\calE(\eta)$.
For example $S(\vphi\bet\gam)=\{\bet,\gam\}$.
$S(0):=S(\mK):=\emptyset$ and $S(\eta)=\{\eta\}$ when $\eta\in\calE(\eta)$.}
\edf

\bprp\label{prp:G}
For $\alp,\kap,a,b\in OT$,
\benu
\item\label{prp:G1}
$G_{\kap}(\alp)\leq\alp$.

\item\label{prp:G2}
$\alp\in\calh_{a}(b) \Rarw G_{\kap}(\alp)\subset\calh_{a}(b)$.

\item\label{prp:GF}
Let $\gam\leq\del$.
Then
$F_{\gam}(\alp)<\bet \spand F_{\del}(\alp)<\gam \Rarw F_{\del}(\alp)<\bet$.
\eenu
\eprp
{\bf Proof} by simultaneous induction on $\ell\alp$.
It is easy to see that
\beqn\label{eq:G}
G_{\kap}(\alp)\ni\bet \Rarw \bet\prec\kap \spand \ell\kap<\ell\bet\leq \ell\alp
\eeqn
\ref{prp:G}.\ref{prp:G1}.
Consider the case $\alp=\psi_{\pi}^{\vec{\nu}}(a)$ with $\pi\not\preceq\kap$.
First let $\kap<\pi$. Then $G_{\kap}(\alp)=G_{\kap}(\{\pi,a\}\cup K(\vec{\nu}))$.
On the other hand we have  $\fal \gam\in K(\vec{\nu})\cup\{\pi,a\}(K_{\alp}(\gam)<a)$, i.e,
$K(\vec{\nu})\cup\{\pi,a\}\subset\calh_{a}(\alp)$.
Proposition \ref{prp:G}.\ref{prp:G2} with (\ref{eq:G}) 
yields $G_{\kap}(K(\vec{\nu})\cup\{\pi,a\})\subset\calh_{a}(\alp)\cap\kap\subset \calh_{a}(\alp)\cap\pi\subset\alp$.
Hence $G_{\kap}(\alp)<\alp$.

Next let $\pi<\kap$ and $\pi\not\preceq\kap$. Then $G_{\kap}(\alp)=G_{\kap}(\pi)$.
By IH we have $G_{\kap}(\pi)\leq\pi$, and $G_{\kap}(\pi)<\pi$ by $\pi\not\preceq\kap$.
On the other hand we have  $K_{\alp}(\pi)<a$, i.e,
$\pi\in\calh_{a}(\alp)$.
Proposition \ref{prp:G}.\ref{prp:G2} 
yields $G_{\kap}(\pi)\subset\calh_{a}(\alp)\cap\pi\subset\alp$.
Hence $G_{\kap}(\alp)<\alp$.
\\
\ref{prp:G}.\ref{prp:G2}.
Since $G_{\kap}(\alp)\leq\alp$ by Proposition \ref{prp:G}.\ref{prp:G1}, we can assume $\alp\geq b$.
Again consider the case $\alp=\psi_{\pi}^{\vec{\nu}}(a)$ with $\pi\not\preceq\kap$.
Then $K(\vec{\nu})\cup\{\pi,a\}\subset\calh_{a}(b)$ and $G_{\kap}(\alp)\subset G_{\kap}(K(\vec{\nu})\cup\{\pi,a\})$.
IH yields the lemma.
\\
\ref{prp:G}.\ref{prp:GF}.
This is seen by induction on $\ell\alp$.
\eprf

\bprp\label{lem:5uv.1-1}
Let $\bet\preceq\alp=\psi_{\pi}^{\vec{\nu}}(a)$.
Then $F_{\pi}(K(\vec{\nu}))<\bet$.
\eprp
\bprf
Let $pd^{(i-1)}(\bet)=\pi_{i-1}=\psi_{\pi_{i}}^{\vec{\nu}_{i}}(a_{i})$
with $\bet=\pi_{0}$ and $\pi=\pi_{n}$.
Then by $\pi_{i-1}<\pi_{i}$ we have $\pi_{i}\in\calh_{a_{j+1}}(\pi_{j})$ for any $j<i$, and
$K(\vec{\nu})\subset\calh_{a_{j+1}}(\pi_{j})$ 
for $\vec{\nu}=\vec{\nu}_{n}$ and any $j<n$.
On the other hand we have $\calh_{a_{j+1}}(\pi_{j})\cap\pi_{j+1}\subset\pi_{j}$.
We see by induction on $n-j\geq 0$ that
$F_{\pi}(K(\vec{\nu}))<\pi_{j}$.
\eprf

\bprp\label{prp:G4}
Let $\gam\preceq\tau$ and $\gam\not\prec\kap$. Then $G_{\kap}(\tau)\subset G_{\kap}(\gam)$.
\eprp
\bprf
Let  $\gam\not\prec\kap$.
We show $\gam\preceq\tau \Rarw G_{\kap}(\tau)\subset G_{\kap}(\gam)$
by induction on $\ell\gam-\ell\tau$.
Let $\gam\preceq\tau=\psi_{\pi}^{\vec{\nu}}(a)$.
By IH we have $G_{\kap}(\tau)\subset G_{\kap}(\gam)$.
On the other hand we have $G_{\kap}(\pi)\subset G_{\kap}(\tau)$
since $\pi\not\prec\kap$ and $\pi=\kap \Rarw G_{\kap}(\pi)=\emptyset$, cf. (\ref{eq:G}).
\eprf

\bprp\label{prp:G3}
Let $a,\alp,\kap,\bet,\del\in OT$ with $\alp=\psi_{\pi}^{\vec{\nu}}(a)$ for some $\{a\}\cup K(\vec{\nu})\subset OT$.
If
$\bet\not\in\calh_{a}(\alp)$ and $K_{\del}(\bet)<a$, then there exists a $\gam\in F_{\del}(\bet)$
such that $\calh_{a}(\alp)\not\ni\gam<\del$.
\eprp
\bprf
By induction on $\ell\bet$. Assume $\bet\not\in\calh_{a}(\alp)$ and $K_{\del}(\bet)<a$.
By IH we can assume that $\bet=\psi_{\kap}^{\vec{\xi}}(b)$.
If $\bet<\del$, then $\bet\in F_{\del}(\bet)$, and $\gam=\bet$ is a desired one.
Assume $\bet\geq\del$. Then we have $K_{\del}(\bet)=\{b\}\cup K_{\del}(\{b,\kap\}\cup K(\vec{\xi}))<a$.
In particular $b<a$, and hence $\{b,\kap\}\cup K(\vec{\xi})\not\subset\calh_{a}(\alp)$.
By IH there exists a $\gam\in F_{\del}(\{b,\kap\}\cup K(\vec{\xi}))=F_{\del}(\bet)$
such that $\calh_{a}(\alp)\not\ni\gam<\del$.
\eprf


\section{Distinguished sets}\label{sec:distinguished}

In this section, working in the set theory KP$\ell$ for limits of admissibles,
we will develop rudiments of distinguished classes, which was first introduced by W. Buchholz\cite{Buchholz75}.
Since many properties of distingusihed classes are seen as in \cite{Wienpi3d, WFnonmon2}, we will omit their proofs.

As in \cite{WFnonmon2} our welfoundedness proof inside KP$\Pi_{N}$ goes as follows.
The wellfoundedness of $OT$ is reduced to one of the relation $\prec$ in the following way.
$\alp\in V(X)$ in Definition \ref{df:CX}.\ref{df:CXV} is intended for $\alp$ 
to be in the wellfounded part of $\prec$ with respect to a set $X$.
In Lemma \ref{th:3wf16} it is shown for a $\Del_{1}$ class $\calg(X)$ defined in Definition \ref{df:calg},
that $\eta\in\calg(X)\cap V(X)$ yields the existence of a distinguished set $X^{\prime}$
such that $\eta\in X^{\prime}$ \textit{provided that} $X$ is a distinguished set which is closed under the `hyperjump' operation 
$X\mapsto X^{\prime}$ for any $\gam\prec\eta$.
Let us call such an $X$ \textit{$\eta$-Mahlo}.
It turns out that we need the fact that $X\subset V(X)$ for any distinguished sets $X$ in proving Lemma \ref{th:3wf16}.
Furthermore we need even stronger condition $X\subset V^{*}(X)$ for the Claim \ref{clm:3wf1632} in Lemma \ref{th:3wf16},
where $V^{*}(X)$ is defined in Definition \ref{df:CX}.\ref{df:CXV*}.
This motivates our Definition \ref{df:3wfdtg32}.\ref{df:3wfdtg.832} of distinguished sets (\ref{eq:distinguishedclass}).

There remain three tasks for each $\eta\in OT$.
One is to show that $\eta\in\calg(X)$, second to show $\eta\in V(X)$, and third the existence of an $\eta$-Mahlo
distinguished set.
It is not hard to show $\eta\in\calg(X)$ by induction on $a$ for $\eta=\psi_{\pi}^{\vec{\nu}}(a)$, cf. Lemma \ref{th:id5wf21}.
Next for sets $P$ let $\calw^{P}$ be the maximal distinguished class \textit{in} $P$.
$\calw^{P}$ is $\Sig_{1}^{P}$, i.e., $\Sig_{1}$-definable class on $P$, and $\calw^{Q}$ is a distinguished \textit{set} in $P$
for any sets $Q\in P$, cf. subsection \ref{subsec:Mahlouniv}.
In particular $\calw=\calw^{\rm{L}}$ is the maximal distinguished class for the whole $\Pi_{N}$-reflecting universe L.
Let us say that $P$ is $\eta$-Mahlo if $\calw^{P}$ is an $\eta$-Mahlo distinguished class.
In view of Lemma \ref{th:3wf16} $P$ is $\eta$-Mahlo if $P$ is $\Pi_{2}$-reflecting on $\gam$-Mahlo sets for any 
$\gam\prec\eta$ since $\calg(\calw^{P})$ is $\Pi_{2}^{P}$.
This means that we need to iterate recursively Mahlo operations along $\prec$ up to a given $\eta$
assuming that $\eta$ is in the wellfounded part $V(\calw)$.
Now if $\gam\prec\eta$, then the sequence of ordinals $\{m_{k}(\gam)\}_{k}$ is smaller than $\{m_{k}(\eta)\}_{k}$
in a sense. Indeed we could assign an ordinal $o_{1}(\{m_{k}(\gam)\})<\veps_{\mK+2}$
in such a way that $o_{1}(\{m_{k}(\gam)\})<o_{1}(\{m_{k}(\eta)\})$ as in Definition \ref{df:lxo}.
However if we refer such a big ordinal $o_{1}(\{m_{k}(\eta)\})>\eta$ explicitly in
defining $\eta$ to be in $V(\calw)$,
the persistency (\ref{eq:Vpersistency}) in Definition \ref{df:CX}.\ref{df:CXV} would be lost.
As we see it in this section, the persistency is crucial for distinguished sets, cf. Proposition \ref{lem:3wf5}.

The $k$-predecessors are needed for us to embed the relation $\prec$ on $OT$ to an exponential structure
induced solely from ordinals $\{m_{k}(\alp)\}_{k}$ (cf. Lemma \ref{lem:exponentialinprec}), which in turn yields sets $V(X)=V_{N}(X)$ introduced in subsection \ref{subsec:V(X)}
with the persistency (\ref{eq:Vpersistency}) in Definition \ref{df:CX}.


 $X,Y,\ldots$ range over {\it subsets\/} of $OT_{n}$. 
 While $\calx,\caly,\ldots$ range over {\em classes}.

We define sets $\calc^{\alp}(X)\subset OT_{n}$ for $\alp\in OT_{n}, X\subset OT_{n}$ as follows.

\bdf\label{df:CX}
{\rm Let} $\alp,\bet\in OT_{n}, X\subset OT_{n}$.
\benu
\item
{\rm Let} 
$\calc^{\alp}(X)$ {\rm be the closure of} $\{0,\mK\}\cup(X\cap\alp)$
{\rm under} $+$, $\mK<\bet\mapsto\ome^{\bet}\in OT_{n}$,
$(\bet,\gam)\mapsto\vphi\bet\gam\,(\bet,\gam<\mK)$,
$\mK> \bet\mapsto\Ome_{\bet}>\bet$,
{\rm and} $(\sig,\bet,\vec{\xi})\mapsto \psi_{\sig}^{\vec{\xi}}(\bet)$
 {\rm for} $ \sig>\alp$ {\rm in}  $OT_{n}$.

{\rm The last clauses say that, if}
$\Ome_{\bet}>\bet\in\calc^{\alp}(X) \Rarw \Ome_{\bet}\in\calc^{\alp}(X)$,
{\rm and}
$\psi_{\sig}^{\vec{\xi}}(a)\in\calc^{\alp}(X)$ {\rm if}
$\{\sig,a\}\cup K(\vec{\xi})\subset\calc^{\alp}(X)$ {\rm and} $\sig>\alp$.

\item
$\alp^{+}=\Ome_{a+1}$ {\rm denotes the least recursively regular term above} $\alp$ {\rm if such a term exists. Otherwise} $\alp^{+}:=\infty$.
{\rm Obviously $\alp^{+}$ is computable from $\alp$.}

\item\label{df:CXV}
$V(X)$ {\rm is a $\Del_{1}$-class such that}
\beqnarr
&&
\fal\alp<\mK[(X\cap\alp=Y\cap\alp \Rarw V(X)\cap\alp^{+}=V(Y)\cap\alp^{+}) 
\label{eq:Vpersistency}
\\
& \land & (\lnot\exi\kap,a,\vec{\xi}\neq\vec{0}(\alp=_{NF}\psi_{\kap}^{\vec{\xi}}(a)) \Rarw \alp\in V(X))]
\nonumber
\eeqnarr

\item
$V\calc^{\alp}(X):=V(X)\cap\calc^{\alp}(X)$.

\item\label{df:CXV*}
$\alp\in V^{*}(X) :\Lrarw \alp\in V(X) \spand \calc^{\alp}(X)\cap\alp\subset V(X)$.

\item
$V^{*}\calc^{\alp}(X):=V^{*}(X)\cap\calc^{\alp}(X)$.

\eenu
\edf

\bprp\label{lem:CX1}
$X\cap\alp=Y\cap\alp \Rarw \calc^{\alp}(X)=\calc^{\alp}(Y)$ and $X\mapsto\calc^{\alp}(X)$ is monotonic.
\eprp

\bprp\label{prp:CX0}
$\alp<\bet<\alp^{+} \Rarw \calc^{\alp}(X)\subset\calc^{\bet}(X)$.
\eprp
\bprf
By induction on $\ell\gam\,(\gam\in OT_{n})$ we see that
$\gam\in\calc^{\alp}(X) \Rarw \gam\in\calc^{\bet}(X)$.
\eprf

\bprp\label{prp:kFC}
Let $\del\leq\mK$. Then
$F_{\del}(\alp)\cup k_{\del}(\alp)\subset X \Rarw \alp\in\calc^{\mK}(X)$.
\eprp
\bprf
This is seen by induction on $\ell\alp$.
\eprf

\bprp\label{lem:CX4}
Assume $\alp\in\calc^{\alp}(X)$ and $\alp\preceq\sig$.
Then $\sig\in\calc^{\alp}(X)$.
\eprp
\bprf
We see by induction on $\ell\alp-\ell\sig$ that
$\alp\in\calc^{\alp}(X) \spand \alp\preceq\sig \Rarw \sig\in\calc^{\alp}(X)$.
\eprf

\bprp\label{lem:CX2} {\rm (Cf. \cite{Wienpi3d}, Lemmas 3.5.3 and 3.5.4.)}
Assume $\fal\gam\in X[ \gam\in\calc^{\gam}(X)]$ for a set $X\subset OT_{n}$.

\benu
\item\label{lem:CX2.3} 
$\alp\leq\bet \Rarw \calc^{\bet}(X)\subset \calc^{\alp}(X)$.

\item\label{lem:CX2.4} 
$\alp<\bet<\alp^{+} \Rarw \calc^{\bet}(X)=\calc^{\alp}(X)$.
\eenu
\eprp
\bprf
\ref{lem:CX2}.\ref{lem:CX2.3}.
We see by induction on $\ell\gam\,(\gam\in OT_{n})$ that
\beqn\label{eq:CX2.3}
\fal\bet\geq\alp[\gam\in \calc^{\bet}(X) \Rarw \gam\in \calc^{\alp}(X)\cup(X\cap\bet)]
\eeqn
For example, if $\psi_{\pi}^{\vec{\xi}}(\del)\in \calc^{\bet}(X)$ with $\pi>\bet\geq\alp$
and $\{\pi,\del\}\cup K(\vec{\xi})\subset\calc^{\alp}(X)\cup(X\cap\bet)$, then 
$\pi\in\calc^{\alp}(X)$, and 
for any $\gam\in\{\del\}\cup K(\vec{\xi})$, either $\gam\in\calc^{\alp}(X)$ or
$\gam\in X\cap\bet$ by IH. If $\gam<\alp$, then $\gam\in X\cap\alp\subset\calc^{\alp}(X)$.
If $\alp\leq\gam\in X\cap\bet$, then $\gam\in\calc^{\gam}(X)$ by the assumption, and
by IH we have $\gam\in\calc^{\alp}(X)\cup(X\cap\gam)$, i.e., $\gam\in\calc^{\alp}(X)$.
Therefore $\{\pi,\del\}\cup K(\vec{\xi})\subset\calc^{\alp}(X)$, and 
$\psi_{\pi}^{\vec{\xi}}(\del)\in\calc^{\alp}(X)$.

Using (\ref{eq:CX2.3}) we see from the assumption that
$
\fal\bet\geq\alp[ \gam\in\calc^{\bet}(X) \Rarw \gam\in\calc^{\alp}(X)]
$.
\\

\noindent
\ref{lem:CX2}.\ref{lem:CX2.4}.
Assume $\alp<\bet<\alp^{+}$. Then by Proposition \ref{lem:CX2}.\ref{lem:CX2.3} we have
$\calc^{\bet}(X)\subset\calc^{\alp}(X)$.
Conversely $\calc^{\alp}(X) \subset\calc^{\bet}(X)$ is seen from Proposition \ref{prp:CX0}.
\eprf

\bdf\label{df:wftg}
\benu
\item $Prg[X,Y] :\Lrarw \fal\alp\in X(X\cap\alp\subset Y \to \alp\in Y)$.

\item {\rm For a definable class} $\calx$, $TI[\calx]$ {\rm denotes the schema:}\\
$TI[\calx] :\Lrarw Prg[\calx,\caly]\to \calx\subset\caly \mbox{ {\rm holds for} any definable class } \caly$.
\item
{\rm For} $X\subset OT_{n}$, $W(X)$ {\rm denotes the} wellfounded part {\rm of} $X$. 
\item 
$Wo[X] : \Lrarw X\subset W(X)$.
\eenu
\edf
Note that for $\alp\in OT_{n}$,
$W(X)\cap\alp=W(X\cap\alp)$.

\bdf \label{df:3wfdtg32}
{\rm For} $X\subset OT_{n}$ {\rm and} 
$\alp\in OT_{n}$,
\benu
\item\label{df:3wfdtg.832}
\beqn\label{eq:distinguishedclass}
D[X] :\Lrarw X<\mK \spand
\fal\alp(\alp\leq X\to W(V^{*}\calc^{\alp}(X))\cap\alp^+= X\cap\alp^+)
\eeqn

{\rm A class} $\calx$ {\rm is said to be a} distinguished class {\rm if} $D[\calx]${\rm . A} distinguished set {\rm is a set which is a distinguished class.}

\item\label{df:3wfdtg.932}
$\calw:=\bigcup\{X :D[X]\}$.


\eenu
\edf

Since, in KP$\ell$,
the wellfounded part $W(X)$ of a set $X$ is again a set, $D[X]$ is
 $\Del_{1}$.
Hence both $\calw$ and $\calc^{\alp}(\calw)$ are $\Sig_{1}$.
Obviously any distinguished set $X$ enjoys the condition $\fal\alp\in X[\alp\in V^{*}\calc^{\alp}(X)]$.

\bprp\label{lem:3wf4}
$D[X]\Rarw Wo[X]$.
\eprp

\bprp\label{lem:3.11.632}{\rm (Cf. Lemma 3.30 in \cite{WFnonmon2}.)}\\
Let $X$ be a distinguished set. Then
$\alp\in X \Rarw \fal\bet[\alp\in\calc^{\bet}(X)]$.
\eprp

\bprp\label{lem:5uv.232general}{\rm (Cf. Lemma 3.28 in \cite{WFnonmon2}.)}\\
For any distinguished sets $X$ and $Y$, the following holds:
\[
X\cap\alp=Y\cap\alp \Rarw \fal\bet<\alp^{+}\{V^{*}\calc^{\bet}(X)\cap\bet^{+}=V^{*}\calc^{\bet}(Y)\cap\bet^{+}\}.
\]
\eprp
\bprf
Assume that $X\cap\alp=Y\cap\alp$ and $\bet<\alp^{+}$.
By the condition (\ref{eq:Vpersistency}) we have
 $V(X)\cap\bet^{+}=V(Y)\cap\bet^{+}$. 

On the other hand we have 
by Propositions \ref{lem:CX2}.\ref{lem:CX2.4} and \ref{lem:CX1}, 
$\calc^{\bet}(X)=\calc^{\bet}(Y)$, and for any $\del<\bet^{+}$,
$\calc^{\del}(X)=\calc^{\del}(Y)$.
Hence  $V^{*}(X)\cap\bet^{+}=V^{*}(Y)\cap\bet^{+}$. 
\eprf

\bprp\label{lem:3wf5}
Let $X$ and $Y$ be distinguished sets.

\benu
\item\label{lem:3wf5.1}
$\alp\leq X \spand \alp\leq Y \Rarw X\cap\alp^+=Y\cap\alp^+$.

\item\label{lem:3wf5.2} 
Either $X\subset_{e} Y$ or $Y\subset_{e} X$, where
$X\subset_{e} Y$ designates that $Y$ is an end extension of $X$, i.e.,
$
X\subset_{e} Y : \Lrarw  X\subset Y \spand \fal\alp\in Y\fal\bet\in X(\alp<\bet\to \alp\in X)
$.
\eenu
\eprp

\bprp\label{lem:3wf6}
$\calw$ is the maximal distinguished class, i.e.,
$D[\calw]$.
Also $TI[\calw]$ for $\calw\subset\mK$.
\eprp

\subsection{Sets $\calc^{\alp}(X)$ and $\calg(X)$}\label{subsec:C(X)}

In this subsection we will establish elementary properties on sets $\calc^{\alp}(X)$.

\bprp\label{lem:CX3} {\rm (Cf. \cite{Wienpi3d}, Lemma 3.6.)}
Let $\gam<\bet$.
For a distinguished set $X$ assume $\alp\in\calc^{\gam}(X)$
 and $\fal\kap\leq\bet[G_{\kap}(\alp)<\gam]$.
\benu
\item\label{lem:CX3.1}
Assume $\mbox{{\rm LIH}} :\fal\del[\ell\del\leq\ell\alp\spand\del\in\calc^{\gam}(X)\cap\gam\Rarw\del\in \calc^{\bet}(X)]$.
Then $\alp\in\calc^{\bet}(X)$.
\item\label{lem:CX3.2}
$\calc^{\gam}(X)\cap\gam \subset X \Rarw \alp\in\calc^{\bet}(X)$.
\eenu
\eprp
\bprf
\ref{lem:CX3}.\ref{lem:CX3.1} by induction on $\ell\alp$. 
If $\alp<\gam$, then $\alp\in \calc^{\gam}(X)\cap\gam$.
LIH yields $\alp\in \calc^{\bet}(X)$. 
Assume $\alp\geq\gam$. 
Except the case $\alp=\psi_{\pi}^{\vec{\nu}}(a)$ for some $\pi,a,\vec{\nu}$, IH yields $\alp\in\calc^{\bet}(X)$. 
Suppose $\alp=\psi_{\pi}^{\vec{\nu}}(a)$ for some $\{\pi,a\}\cup K(\vec{\nu})\subset\calc^{\gam}(X)$ and $\pi>\gam$.
If $\pi\leq\bet$, then $\{\alp\}=G_{\pi}(\alp)<\gam$ by the second assumption. Hence this is not the case, and we obtain 
$\pi>\bet$.
 Then $G_{\kap}(\{\pi,a\}\cup K(\vec{\nu}))=G_{\kap}(\alp)<\gam$ for any $\kap\leq\bet<\pi$. 
IH yields $\{\pi,a\}\cup K(\vec{\nu})\subset\calc^{\bet}(X)$. We conclude $\alp\in\calc^{\bet}(X)$ from $\pi>\bet$.

\eprf

\bdf\label{df:calg}
$\calg(X):=\{\alp:\alp\in \calc^{\alp}(X) \spand \calc^{\alp}(X)\cap\alp\subset X\}$.
\edf

\bprp\label{lem:CMsmpl}
Let $\alp\in\calc^{\bet}(X)$ and $X\cap\bet\subset\calg(X)$ for a distinguished set $X$.
Assume $X\cap\bet<\del$.
Then $F_{\del}(\alp)\subset\calc^{\bet}(X)$.
\eprp
\bprf
By induction on $\ell\alp$.
Let $\{0,\mK\}\not\ni\alp\in\calc^{\bet}(X)$.
First consider the case $\alp\not\in\calE(\alp)$.
If $\alp\in X\cap\bet\subset\calg(X)$, then 
$\calE(\alp)\subset\calc^{\alp}(X)\cap\alp\subset X\subset\calc^{\bet}(X)$ by Proposition \ref{lem:3.11.632}.
Otherwise we have $\alp\not\in\calE(\alp)\subset\calc^{\bet}(X)$.
In each case IH yields $F_{\del}(\alp)=F_{\del}(\calE(\alp))\subset\calc^{\bet}(X)$.

Let $\alp=\psi_{\pi}^{\vec{\nu}}(a)$ for some $\pi,\vec{\nu},a$. 
If $\alp<\del$, then $F_{\del}(\alp)=\{\alp\}$, and there is nothing to prove.
Let $\alp\geq\del$. Then $F_{\del}(\alp)=F_{\del}(\{\pi,a\}\cup K(\vec{\nu}))$.
On the other side we see $\{\pi,a\}\cup K(\vec{\nu})\subset\calc^{\bet}(X)$ from $\alp\in\calc^{\bet}(X)$ 
and the assumption.
IH yields $F_{\del}(\alp)\subset\calc^{\bet}(X)$.
\eprf
\\

Next we show $X\subset\calg(X)$ for any distinguished set $X$, 
cf. Lemma \ref{lem:wf5.332}.

\bprp\label{lem:id7wf8}
Let $X$ be a distinguished set, and assume $X\cap\bet\subset\calg(X)$.
\benu
\item\label{lem:id7wf8.3}
$\fal\tau[\alp\in X\cap\bet \Rarw G_{\tau}(\alp)\subset X]$.

\item\label{lem:id7wf8.3.5}
$\fal\bet\fal\tau[\alp\in\calc^{\bet}(X) \Rarw G_{\tau}(\alp)\subset \calc^{\bet}(X)]$.
\eenu
\eprp
\bprf
 By simultaneous induction on $\ell\alp$.
\\
\ref{lem:id7wf8}.\ref{lem:id7wf8.3}. 
Suppose $\alp\in X\cap\bet\subset\calg(X)$. Then  $\alp\in \calc^{\alp}(X)$, and $\calc^{\alp}(X)\cap\alp\subset X$.

Let $\alp\not\in\calE(\alp)$. Then $\calE(\alp)\subset \calc^{\alp}(X)\cap\alp\subset X$.
 IH yields
$G_{\tau}(\alp)=G_{\tau}(\calE(\alp))\subset X$.
Assume $\alp\in\calE(\alp)$, i.e., $\alp=\psi_{\pi}^{\vec{\nu}}(a)$ for some $\pi,a,\vec{\nu}$.
Then $\{\pi,a\}\cup K(\vec{\nu})\subset\calc^{\alp}(X)$ by $\alp\in\calc^{\alp}(X)$.
We can assume $\pi\not\preceq\tau$. Then $G_{\tau}(\alp)\subset G_{\tau}(\{\pi,a\}\cup K(\vec{\nu}))$.
By IH with Proposition \ref{prp:G}.\ref{prp:G1} we have 
$G_{\tau}(\alp)\subset G_{\tau}(\{\pi,a\}\cup K(\vec{\nu}))\subset\calc^{\alp}(X)\cap\alp\subset X$.
\\

\noindent
\ref{lem:id7wf8}.\ref{lem:id7wf8.3.5}. 
Assume $\alp\in\calc^{\bet}(X)$. 
We show $G_{\tau}(\alp)\subset\calc^{\bet}(X)$. 
If $\alp\in X\cap\bet$, then by Proposition\ref{lem:id7wf8}.\ref{lem:id7wf8.3} we have 
$G_{\tau}(\alp)\subset X\cap\bet\subset \calc^{\bet}(X)$. 
Consider the case $\alp\not\in X\cap\bet$.
If $\alp\not\in\calE(\alp)$, then IH yields $G_{\tau}(\alp)=G_{\tau}(\calE(\alp))\subset\calc^{\bet}(X)$.
Let $\alp=\psi_{\pi}^{\vec{\nu}}(a)$ for some $\{\pi,a\}\cup K(\vec{\nu})\subset\calc^{\bet}(X)$ with $\bet<\pi\not\preceq\tau$.
IH yields $G_{\tau}(\alp)\subset G_{\tau}(\{\pi,a\}\cup K(\vec{\nu}))\subset\calc^{\bet}(X)$.
\eprf

\bprp\label{lem:KC} 
Let $X$ be a distinguished set, and assume $X\cap\bet\subset\calg(X)$. Then
\[
\fal\alp\fal\sig\leq\bet[\alp\in\calc^{\bet}(X)\Rarw G_{\sig}(\alp)\subset X]
.\]
\eprp
\bprf 
By induction on $\ell\alp$ using Proposition \ref{lem:id7wf8}.\ref{lem:id7wf8.3}
 we see $\alp\in\calc^{\bet}(X)\spand\sig\leq\bet \Rarw G_{\sig}(\alp)\subset X$. 
\eprf

\bprp\label{lem:wf5.3.3-132}
Let $X$ be a distinguished set. 
Assume 
$X\cap\gam\subset\calg(X)$, and $\alp\in\calc^{\gam}(X)\cap\gam$.
Then $\calc^{\alp}(X)\cap\alp\subset \calc^{\gam}(X)$.
\eprp
\bprf
First suppose that there exists a $\del$ such that $\alp\leq\del\in X\cap\gam\subset\calg(X)$.
Then $\calc^{\del}(X)\cap\del\subset X$.
If $\alp=\del$, then $\calc^{\alp}(X)\cap\alp\subset X\subset\calc^{\gam}(X)$ by Proposition \ref{lem:3.11.632}.
Let $\alp<\del$.
Then $X\cap\del\subset\calg(X)$, and $\alp\in\calc^{\del}(X)\cap\del$ by Proposition \ref{lem:CX2}.\ref{lem:CX2.3}.
Moreover we have $\del\in X$.
Therefore it suffices to show the proposition under the \textit{assumption} $\gam\in X$,
for then $\calc^{\alp}(X)\cap\alp\subset\calc^{\del}(X)\cap\alp\subset X\subset\calc^{\gam}(X)$.

Let us prove the proposition by main induction on $\gam\in X$.
If $\alp\leq X\cap\gam$, then MIH yields the proposition as we saw it above.
In what follows assume $X\cap\gam<\alp$.

By subsidiary induction on $\ell\alp+\ell\bet$ we show that
\[
\bet\in\calc^{\alp}(X)\cap\alp \Rarw \bet\in\calc^{\gam}(X).
\]
If $\bet\in X$, then $\bet\in\calc^{\gam}(X)$ follows from Proposition \ref{lem:3.11.632}.
In what follows suppose $\bet\not\in X$

If $\bet\not\in\calE(\bet)$, then $\bet\in\calc^{\gam}(X)$ is seen from SIH.
Assume $\bet=\psi_{\pi}^{\vec{\nu}}(a)$ with a $\pi>\alp$ and some $K(\vec{\nu})\cup\{\pi,a\}\subset \calc^{\alp}(X)$.
If $\alp\not\in\calE(\alp)$, then $\bet\leq\del$ for some $\del\in\calE(\alp)\subset\calc^{\gam}(X)\cap\gam$.
Since $\ell\del<\ell\alp$, SIH yields $\bet\in\calc^{\gam}(X)$.
Let $\alp=\psi_{\kap}^{\vec{\xi}}(b)$ for some $\kap,b,\vec{\xi}$.
By $X\not\ni\alp\in\calc^{\gam}(X)$ we have $\gam<\kap$.

First consider the case $\gam<\pi$.
Then $\fal\sig\leq\gam[G_{\sig}(\{\pi,a\}\cup K(\vec{\nu}))=G_{\sig}(\bet)<\bet<\alp]$ by Proposition \ref{prp:G}.\ref{prp:G1}.
Since $\ell\eta<\ell\bet$ for each $\eta\in \{\pi,a\}\cup K(\vec{\nu})$,
by SIH we have LIH: $\fal\del[\ell\del\leq\ell\eta \spand\del\in\calc^{\alp}(X)\cap\alp\Rarw\del\in \calc^{\gam}(X)]$ 
in  Proposition \ref{lem:CX3}.\ref{lem:CX3.1}, which yields $\{\pi,a\}\cup K(\vec{\nu})\subset\calc^{\gam}(X)$,
and $\bet\in\calc^{\gam}(X)$. 

Next assume $\pi\leq\gam<\kap$.
$\pi\not\in\calh_{b}(\alp)$ since otherwise by $\pi<\kap$ we would have $\pi<\alp$.
Then by Proposition \ref{prp:psicomparison} we have $a\geq b$ and $K(\vec{\xi})\cup\{\kap,b\}\not\subset\calh_{a}(\bet)$.
On the other hand we have $K_{\alp}(K(\vec{\xi})\cup\{\kap,b\})<b\leq a$.
By Proposition \ref{prp:G3} pick a $\del\in F_{\alp}(K(\vec{\xi})\cup\{\kap,b\})$
such that $\calh_{a}(\bet)\not\ni\del<\alp$.
We have $\ell\del<\ell\alp$ and $K(\vec{\xi})\cup\{\kap,b\}\subset\calc^{\gam}(X)$.
Hence by Proposition \ref{lem:CMsmpl} we obtain $\del\in\calc^{\gam}(X)\cap\gam$.
From $\bet\not\in\calh_{a}(\bet)$ we see $\bet\leq\del$.
If $\bet=\del\in\calc^{\gam}(X)$, we are done.
Let $\bet<\del$.
Then $\bet\in\calc^{\del}(X)\cap\del$, and SIH with $\ell\del<\ell\alp$ yields $\bet\in\calc^{\gam}(X)$.
\eprf

\blem\label{lem:wf5.332}
Let $X$ be a distinguished set. 
Then $X\subset\calg(X)$, 
$\fal\alp\in X \fal\tau(G_{\tau}(\alp)\subset X)$, and
$\fal\alp\in\calc^{\gam}(X)\cap\gam(\calc^{\alp}(X)\cap\alp\subset \calc^{\gam}(X))$.
\elem
\bprf 
We have $\gam\in\calc^{\gam}(X)$ for $\gam\in X$.
%
%

Assume $\alp\in\calc^{\gam}(X)\cap\gam$. 
We have $\gam\in W(\calc^{\gam}(X))\cap\gam^{+}=X\cap\gam^{+}$ by $\gam\in X$. 
Hence $\alp\in W(\calc^{\gam}(X))\cap\gam^{+}\subset X$.
Next $\fal\alp\in X \fal\tau(G_{\tau}(\alp)\subset X)$ is seen from $X\subset\calg(X)$ and Proposition \ref{lem:id7wf8}.\ref{lem:id7wf8.3}.
Finally $\fal\alp\in\calc^{\gam}(X)\cap\gam(\calc^{\alp}(X)\cap\alp\subset \calc^{\gam}(X))$ is seen from Proposition 
\eprf
\\

The following Propositions \ref{prp:updis} and \ref{prp:maxwup} are seen from Lemma \ref{lem:wf5.332}.

\bprp\label{prp:updis}
Let $X$ be a distinguished set.
Then $\alp\leq X\cap\bet \spand \alp\in\calc^{\bet}(X) \Rarw \alp\in X$.
\eprp

\bprp\label{prp:maxwup}
Let $X$ be a distinguished set, and $\alp\in X$.
Then $\calc^{\alp}(X)\cap \alp\subset X$.
\eprp


\bprp\label{prp:maxwup1}
Let $X$ be a distinguished set, and $\alp\in X$.
Then $S(\alp)\subset X$.
\eprp
\bprf
Let $\alp\in X$.
Then $\alp\in\calc^{\alp}(X)$ by Proposition \ref{lem:3.11.632}.
Hence $S(\alp)\cap\alp\subset\calc^{\alp}(X)\cap\alp\subset X$ by Proposition \ref{prp:maxwup}.
\eprf

\bprp\label{prp:id3wf20a-1}
Let $X$ be a distinguished set.
$\alp\in\calc^{\del}(X)\Rarw F_{\del}(\alp)\subset X$.
\eprp
{\bf Proof} by induction on $\ell\alp$.
If $\alp\in X\cap\del$, then $S(\alp)\subset X$ by Proposition \ref{prp:maxwup1}, and 
$F_{\del}(\alp)=F_{\del}(S(\alp))\subset X$ by IH.
Otherwise $S(\alp)\subset\calc^{\del}(X)$, and $F_{\del}(\alp)=F_{\del}(S(\alp))\subset X$ by IH.
\eprf

\bprp\label{lem:3wf9}
Let $X$ be a distinguished set, and put $Y=W(V^{*}\calc^{\alp}(X))\cap\alp^+$ for an $\alp<\mK$.
Assume that $\alp\in\calg(X)$ and
\[
\fal\bet<\mK(X<\bet \spand \bet^+<\alp^+ \Rarw 
W(V^{*}\calc^{\bet}(X))\cap\bet^{+}\subset X)
.\]
Then $\alp\in Y$ and $D[Y]$.
\eprp
\bprf 
As in \cite{odMahlo, WFnonmon2}  
this is seen from Lemma \ref{lem:wf5.332}.
\eprf


\bprp\label{lem:3wf8}
$0\in X$ for any distinguished set $X\neq\emptyset$.
\eprp
\bprf 
This is seen from Propositions \ref{lem:3wf9} and  \ref{lem:3wf5}.\ref{lem:3wf5.1}.
\eprf
\\

The following Lemma \ref{th:3wf16} is a key on distinguished classes.

\blem\label{th:3wf16}{\rm (Cf. Lemma 3.3.7 in \cite{WFnonmon2}.)}\\
Let $X$ be a distinguished set, and suppose for an $\eta<\mK$
\beqn\label{eq:3wf16hyp.132}
\eta\in\calg(X)\cap V(X)
\eeqn
and 
\beqn\label{eq:3wf16hyp.232}
\fal\gam\prec\eta(\gam\in\calg(X)\cap V(X)  
\to \gam\in X)
\eeqn
 Then
\[
\eta\in W(V^{*}\calc^{\eta}(X))\cap\eta^{+} \mbox{ and } D[W(V^{*}\calc^{\eta}(X))\cap\eta^{+}]
.\]
\elem
\bprf 
By Proposition \ref{lem:3wf9} and the hypothesis (\ref{eq:3wf16hyp.132}) it suffices to show that
\[
\fal\bet<\mK(X<\bet \spand \bet^+<\eta^+ \Rarw 
W(V^{*}\calc^{\bet}(X))\cap\bet^+\subset X)
.\]
Assume $X<\bet<\mK$ and $\bet^+<\eta^+$. 
We have to show  $W(V^{*}\calc^{\bet}(X))\cap\bet^+\subset X$. 
We prove this by induction on $\gam\in W(V^{*}\calc^{\bet}(X))\cap\bet^+$. 
Suppose $\gam\in V^{*}\calc^{\bet}(X)\cap\bet^{+}$ and 
\[
\mbox{{\rm MIH : }} V^{*}\calc^{\bet}(X)\cap\gam\subset X
\]
We show $\gam\in X$. 

First note that
$
\gam\leq X \Rarw \gam\in X
$
since if $\gam\leq \del$ for some $\del\in X$, then by $X<\bet$ and $\gam\in V^{*}\calc^{\bet}(X)$ 
we have $\del<\bet$, 
$\gam\in V^{*}\calc^{\del}(X)$ and $\del\in W(V^{*}\calc^{\del}(X))\cap\del^{+}=X\cap\del^{+}$.
Hence $\gam\in W(V^{*}\calc^{\del}(X))\cap\del^{+}\subset X$.
Therefore we can assume that
\beqn
\label{eq:3wf9hyp.232X}
X<\gam
\eeqn

We show first 
\beqn
\label{eq:3wf9hyp.232}
\gam\in\calg(X)
\eeqn
First $\gam\in \calc^{\gam}(X)$ by $\gam\in \calc^{\bet}(X)\cap\bet^{+}$ and Proposition \ref{lem:CX2}. 
Second we show the following claim by induction on $\ell\alp$:

\bclm\label{clm:3wf1632}
$\alp\in\calc^{\gam}(X)\cap\gam \Rarw  \alp\in X$.
\eclm
{\bf Proof} of Claim \ref{clm:3wf1632}. 
Assume $\alp\in\calc^{\gam}(X)\cap\gam$.
We have $\alp\in V(X)$ by $\gam\in V^{*}(X)$.
Also by Proposition \ref{lem:wf5.332} we have 
$\calc^{\alp}(X)\cap\alp\subset \calc^{\gam}(X)\cap\gam\subset V(X)$.
Hence $\alp\in V^{*}(X)$, and 
We have $\alp\in\calc^{\bet}(X)\cap\gam\Rarw \alp\in X$ by MIH.

We can assume $\gam^{+}\leq\bet$ for otherwise we have 
$\alp\in V^{*}\calc^{\gam}(X)\cap\gam=V^{*}\calc^{\bet}(X)\cap\gam\subset X$ by MIH.
In what follows assume $\alp\not\in X$.

First consider the case $\alp\not\in\calE(\alp)$. 
By induction hypothesis on lengths we have $\calE(\alp)\subset X\subset\calc^{\bet}(X)$, 
and hence $\alp\in V^{*} \calc^{\bet}(X)\cap\gam$.
Therefore $\alp\in X$ by MIH.

In what follows assume $\alp=\psi_{\pi}^{\vec{\nu}}(a)$ for some $\pi>\gam$ such that
 $\{\pi,a\}\cup K(\vec{\nu})\subset\calc^{\gam}(X)$.
\\
{\bf Case 1}. $\bet<\pi$: 
Then $\fal\kap\leq\bet[G_{\kap}(\{\pi,a\}\cup K(\vec{\nu}))=G_{\kap}(\alp)<\alp<\gam]$ by Proposition \ref{prp:G}.\ref{prp:G1}.
 Proposition \ref{lem:CX3}.\ref{lem:CX3.1} with induction hypothesis on lengths yields 
$\{\pi,a\}\cup K(\vec{\nu})\subset\calc^{\bet}(X)$.
Hence $\alp\in V^{*}\calc^{\bet}(X)\cap\gam$ by $\pi>\bet$.
MIH yields $\alp\in X$.
\\

\noindent
{\bf Case 2}. $\bet\geq\pi$: 
We have $\alp<\gam<\pi\leq\bet$. 
It suffices to show that $\alp\leq X$.
Then by (\ref{eq:3wf9hyp.232X}) we have $\alp\leq\del\in X$ for some $\del<\gam$.
$V^{*}\calc^{\del}(X)\ni\alp\leq\del\in X\cap\del^{+}=W(V^{*}\calc^{\del}(X))\cap\del^{+}$
yields $\alp\in W(V^{*}\calc^{\del}(X))\cap\del^{+}\subset X$.

Consider first the case $\gam\not\in\calE(\gam)$.
By Proposition \ref{lem:3wf8} and $\gam<\bet^{+}<\mK$ we can assume that $\gam\not\in\{0,\mK\}$.
Then let $\del=\max S(\gam)$ denote the largest immediate subterm of $\gam$.
Then 
$\del\in\calc^{\gam}(X)\cap\gam\subset V(X)$ by $\gam\in V^{*}\calc^{\gam}(X)$, and 
by (\ref{eq:3wf9hyp.232X}), $X<\gam\in\calc^{\bet}(X)$ we have 
$\del\in \calc^{\bet}(X)\cap\gam$.
Moreover by Lemma \ref{lem:wf5.332} we have $\calc^{\del}(X)\cap\del\subset\calc^{\gam}(X)\cap\gam\subset V(X)$, 
and $\del\in V^{*}\calc^{\bet}(X)\cap\gam$.
Hence $\del\in X$ by MIH.
Also by $\Ome_{\alp}=\alp$, we have $\alp\leq\del$, i.e., $\alp\leq X$, and we are done.

Let $\gam=\psi_{\kap}^{\vec{\xi}}(b)$ for some $b,\vec{\xi}$ and $\kap>\bet$ by (\ref{eq:3wf9hyp.232X}).
We have $\alp<\gam<\pi\leq\bet<\kap$.
$\pi\not\in\calh_{b}(\gam)$ since otherwise by $\pi<\kap$ we would have $\pi<\gam$.
Then by Proposition \ref{prp:psicomparison} we have $a\geq b$ and $K(\vec{\xi})\cup\{\kap,b\}\not\subset\calh_{a}(\alp)$.
On the other hand we have $K_{\gam}(K(\vec{\xi})\cup\{\kap,b\})<b\leq a$.
By Proposition \ref{prp:G3} pick a $\del\in F_{\gam}(K(\vec{\xi})\cup\{\kap,b\})$
such that $\calh_{a}(\alp)\not\ni\del<\gam$.
Also we have $K(\vec{\xi})\cup\{\kap,b\}\subset\calc^{\bet}(X)$.
Hence by Proposition \ref{lem:CMsmpl} we obtain $\del\in\calc^{\bet}(X)\cap\gam$.
Moreover by Lemma \ref{lem:wf5.332} we have $\calc^{\del}(X)\cap\del\subset\calc^{\gam}(X)\cap\gam\subset V(X)$, 
and $\del\in \calc^{\bet}(X)\cap\gam\subset\calc^{\gam}(X)\cap\gam\subset V(X)$.
Hence $\del\in V^{*}\calc^{\bet}(X)\cap\gam$.
Therefore $\alp\leq\del\in X$ by MIH.
We are done.

Thus Claim \ref{clm:3wf1632} is shown.
\eprf
\\

\noindent
Hence we have (\ref{eq:3wf9hyp.232}), $\gam\in\calg(X)\cap V(X)$.
We have $\gam<\bet^{+}\leq\eta \spand \gam\in\calc^{\gam}(X)$.
If $\gam\prec\eta$, then the hypothesis (\ref{eq:3wf16hyp.232}) yields $\gam\in X$.
In what follows assume $\gam\not\prec\eta$.

If $\fal\tau\leq\eta[G_{\tau}(\gam)<\gam]$, then  Proposition \ref{lem:CX3}.\ref{lem:CX3.2} yields 
$\gam\in\calc^{\eta}(X)\cap\eta\subset X$ by $\eta\in\calg(X)$.

Suppose $\exi\tau\leq\eta[G_{\tau}(\gam)=\{\gam\}]$. 
This means, by $\gam\not\prec\eta$, that
$\gam\prec\tau$ for a $\tau<\eta$.
Let $\tau$ denote the maximal such one.
We have $\gam<\tau<\eta$.
Proposition \ref{lem:CX4} with $\gam\in\calc^{\gam}(X)$ yields $\tau\in\calc^{\gam}(X)$.

Next we show that
\beqn\label{eq:3wf1632last}
\fal\kap\leq\eta[G_{\kap}(\tau)<\gam]
\eeqn
Let $\kap\leq\eta$. 
If $\gam\not\prec\kap$, then $G_{\kap}(\tau)\subset G_{\kap}(\gam)<\gam$ by 
Propositions \ref{prp:G4} and \ref{prp:G}.\ref{prp:G1}.
If $\gam\prec\kap$, then by the maximality of $\tau$ we have $\kap\preceq\tau$, and hence
$G_{\kap}(\tau)=\emptyset$, cf.  (\ref{eq:G}).
(\ref{eq:3wf1632last}) is shown.

Hence  Proposition \ref{lem:CX3}.\ref{lem:CX3.2} yields 
$\tau\in\calc^{\eta}(X)$, and
$\tau\in\calc^{\eta}(X)\cap\eta\subset X$ by $\eta\in\calg(X)$.
Therefore $X<\gam<\tau\in X$.
This is not the case by (\ref{eq:3wf9hyp.232X}).
We are done.
\eprf

\subsection{Mahlo universes}\label{subsec:Mahlouniv}

\bdf\label{df:3auni}
\benu
\item 
{\rm By a} universe {\rm we mean either a} whole universe {\rm L or a transitive set} 
$Q\in \mbox{{\rm L}}$ {\rm in a whole universe L such that} $\ome\in Q${\rm . Universes are denoted} $P,Q,\ldots$

\item 
{\rm A universe} $P$ {\rm is said to be a} limit universe {\rm if} 
$P$ {\rm is a limit of admissible sets.}
$Lmtad$ {\rm denotes the class of limit universes.}


\item 
{\rm For a universe} $P$, $\Del_{0}(\Del_{1})$ in $P$ {\rm denotes the class of predicates which are} $\Del_{0}$ {\rm in some} $\Del_{1}$ {\rm predicates on} $P${\rm .}

\item
$\alp\in rM_{i}(X)   :\Lrarw   \alp \mbox{ {\rm is} } \Pi_{i}\mbox{{\rm -reflecting on} } X$.
\eenu
\edf

We see the absoluteness of the predicate $D[X]$ over limit universes.

\bprp\label{lem:3ahier} 
Let $P$ be a limit universe and $X\in {\cal P}(\ome)\cap P$. 
\benu
\item 
$W(V^{*}\calc^{\alp}(X))$ is $\Del_{1}$ and $D[X]$ is $\Del_{0}(\Del_{1})$.
\label{lem:3ahier.1}
\item 
$W(V^{*}\calc^{\alp}(X))=\{\alp : P\models \alp\in W(V^{*}\calc^{\alp}(X))\}$ 
and $D[X]\Lrarw P\models D[X]$.
\label{lem:3ahier.2}
\eenu
\eprp

\bdf\label{df:3awp}
{\rm For a limit universe} $P$ {\rm set}
\[
\calw^{P}=\bigcup\{X\in P:D[X]\}=\bigcup\{X\in P:P\models D[X]\}
.\]
\edf

Thus $\calw^{\mbox{\footnotesize L}}=\calw$ for the whole universe L.

\bprp\label{lem:ausinP} 
For any limit universe $P$,
$
D[\calw^{P}]
$.
\eprp

\bprp\label{lem:5wuv}
For limit universes $P,Q$,
$
Q\in P\Rarw \calw^{Q}\subset\calw^{P}\spand \calw^{Q}\in P
$.
\eprp

\bprp\label{lem:3afin}
For any limit universe $P$
\[
\bet\in \calc^{\alp}(\calw^{P}) \lrarw 
\exi X\in P\{D[X]\spand \bet\in \calc^{\alp}(X)\}.
\]
\eprp


In the following Proposition \ref{lem:4acalg} by a \textit{$\Pi^{1}_{0}$-class}
we mean a first-order definable class.

\bprp\label{lem:4acalg} 
Let $\calx$ be a $\Pi^{1}_{0}$-class such that $\calx\subset Lmtad$.
Suppose $P\in rM_{2}(\calx)$ and $\alp\in\calg(\calw^{P})$. 
Then
there exists a universe $Q\in P\cap\calx$ such that 
$\alp\in\calg(\calw^{Q})$. 
\eprp
\bprf 
This is seen as in \cite{WFnonmon2}.
\eprf
\\

Lemma \ref{th:3wf16} together with Proposition \ref{lem:4acalg} yields the following 
Corollary \ref{cor:3wf16}, which is the key 
in our wellfoundedness proofs by distinguished sets.

\bcor\label{cor:3wf16}{\rm (Cf. Lemma 6.1 in \cite{LMPS}.)}

Let $\calx$ be a $\Pi^{1}_{0}$-class such that $\calx\subset Lmtad$.
Suppose $P\in rM_{2}(\calx)$ and $\eta\in\calg(\calw^{P})\cap V(\calw^{P})\cap\mK$. 

Assume that there exists a distinguished set $X_{1}\in P$ such that
\beqn\label{eq:Vlocalize}
\fal Q\in P\cap\calx[X_{1}\in Q \Rarw \eta\in V(\calw^{Q})]
\eeqn

Further assume that any $Q\in P\cap\calx$ 
with $X_{1}\in Q$
enjoys the following condition:
\begin{equation}\label{eq:etamahlo}
\forall\gamma\prec\eta\{\gam\in\calg(\calw^{Q})\cap V(\calw^{Q})
 \Rightarrow   \gamma\in \calw^{Q}\}
\end{equation}
Then $\eta\in\calw^{P}$.
\ecor

\bcor\label{cor:3wf16a}
Suppose $L\in rM_{2}(rM_{2}(Lmtad))$ and $S(\eta)\not\ni\eta\in\calg(\calw)\cap\mK$. 
Then $\eta\in\calw$.
\ecor
\bprf
(\ref{eq:Vlocalize}) and $\eta\in V(\calw)$ holds by the condition (\ref{eq:Vpersistency}).
Also any set $Q\in rM_{2}(Lmtad)$ enjoys (\ref{eq:etamahlo}) even if $\eta=_{NF}\Ome_{a+1}$.
Specifically we have $\forall\gamma\prec\eta\{\gam\in\calg(\calw^{Q})\cap V(\calw^{Q})$.
This is seen from Corollary \ref{cor:3wf16} since there is no $\del\prec\gam$.
Hence Corollary \ref{cor:3wf16} yields $\eta\in\calw$.
\eprf

\bprp\label{lem:id3wf19-1}
Suppose $L\in rM_{2}(rM_{2}(Lmtad))$.
Let $\eta<\Lam$ and $S(\eta)\subset \calw$.
Then $\eta\in\calw$.

Specifically
\benu
\item\label{lem:id3wf19-1.1}
$\eta=_{NF}\eta_{m}+\cdots+\eta_{0}\spand \{\eta_{i}:i\leq m\}\subset\calw \Rarw \eta\in\calw\, (m>0)$.

\item\label{lem:id3wf19-1.2}
$\eta=_{NF}\vphi\bet\gam \spand \{\bet,\gam\}\subset\calw \Rarw \eta\in\calw$.

\item\label{lem:id3wf19-1.3}
$\eta=_{NF}\Ome_{a} \spand a\in\calw \Rarw \eta\in\calw$.
\eenu
\eprp
\bprf
We can assume that $\eta\not\in\calE(\eta)$ and $\eta\neq 0$ by Proposition \ref{lem:3wf8}.
We have $S(\eta)\subset\calc^{\eta}(\calw)$ by Proposition \ref{lem:3.11.632}, and hence $\eta\in\calc^{\eta}(\calw)$.
By Corollary \ref{cor:3wf16a} it suffices to show
\beqn\label{eq:id3wf19-1}
\alp\in\calc^{\eta}(\calw)\cap\eta \Rarw \alp\in \calw
\eeqn
\ref{lem:id3wf19-1}.\ref{lem:id3wf19-1.1}.
It suffices to show that
\[
\eta=\bet\dot{+}\gam \spand \{\bet,\gam\}\subset\calw \Rarw \eta\in\calw
\]
by induction on $\gam\in\calw$, where $\bet\dot{+}\gam$ designates the fact that the natural sum 
$\bet\#\gam=\bet+\gam$, and 
$\bet\dot{+}\gam$ denotes the sum $\bet+\gam$.
We have $\eta\in\calc^{\bet}(\calw)=\calc^{\eta}(\calw)$.
We show (\ref{eq:id3wf19-1}).
If $\alp<\bet$, then Proposition \ref{prp:maxwup} yields $\alp\in \calw$.
Let $\alp=\bet\dot{+}\del$ with $\del<\gam$.
Proposition \ref{prp:maxwup} yields $\del\in\calw$. IH yields $\alp\in\calw$.
\\

\noindent
\ref{lem:id3wf19-1}.\ref{lem:id3wf19-1.2}.
By main induction on $\bet\in \calw$ with subsidiary induction on $\gam\in \calw$
we show $\eta=\vphi\bet\gam\in \calw$.
We show (\ref{eq:id3wf19-1}) by induction on $\ell\alp$.
If $\alp=_{NF}\alp_{m}+\cdots+\alp_{0}\, (m>0)$, then the induction hypothesis on the lengths yields
$\{\alp_{i}:i\leq m\}\subset \calw$. 
By Proposition \ref{lem:id3wf19-1}.\ref{lem:id3wf19-1.1} we obtain $\alp\in\calw$.

If $\alp=_{NF}\Ome_{a}$, then $\alp\leq\max\{\bet,\gam\}$.
Proposition \ref{prp:maxwup} yields $\alp\in \calw$.

Finally let $\alp=_{NF}\vphi\bet_{1}\gam_{1}$.
The induction hypothesis on the lengths yields $\{\bet_{1},\gam_{1}\}\subset \calw$.
If $\bet_{1}<\bet$, then MIH yields $\alp\in \calw$.
If $\bet_{1}=\bet$, then $\gam_{1}<\gam$, and SIH yields $\alp\in \calw$.
If $\bet_{1}>\bet$, then $\alp<\gam$.
Proposition \ref{prp:maxwup} yields $\alp\in \calw$.
\\

\noindent
\ref{lem:id3wf19-1}.\ref{lem:id3wf19-1.3}.
By induction on $a\in \calw$ we show $\eta=_{NF}\Ome_{a}\in\calw$.
We show (\ref{eq:id3wf19-1}) by induction on $\ell\alp$.
If either $\alp=_{NF}\alp_{m}+\cdots+\alp_{0}\, (m>0)$ or $\alp=_{NF}\vphi\bet\gam$, then 
the induction hypothesis on the lengths yields $S(\alp)\subset\calw$.
By Propositions \ref{lem:id3wf19-1}.\ref{lem:id3wf19-1.1} and \ref{lem:id3wf19-1}.\ref{lem:id3wf19-1.2}
we obtain $\alp\in\calw$.
Let $\alp=_{NF}\Ome_{b}$. Then $b\in\calw\cap a$, and IH yields $\alp\in\calw$.
\eprf

\section{Iterating recursively Mahlo operations}\label{sect:iterateM}

As in \cite{LMPS, WFnonmon2} we define a tower relation from relations $\{<_{i}: 2\leq i\leq N-1\}$ on $\ome$,

\begin{definition}\label{df:LEe}
{\rm Let} $<_{1}, <_{0}$ {\rm be two transitive relations on} $\omega$.
\begin{enumerate}
\item
{\rm The relation} $<_{E}=E(<_{1},<_{0})$ 
 {\rm  is on sequences} $\langle (n^{1}_{i}, n^{0}_{i}): i<\ell\rangle$ 
{\rm of pairs with} $<_{1}${\rm -decreasing first components} 
 {\rm (}$n^{1}_{i+1}<_{1}n^{1}_{i}${\rm ), and is defined by}
 \begin{eqnarray*}
 && \langle (n^{1}_{i}, n^{0}_{i}): i<\ell_{0}\rangle <_{E}
 \langle (m^{1}_{i}, m^{0}_{i}): i<\ell_{1}\rangle
 \mbox{ {\rm iff}}
 \\
 \mbox{{\rm either}} &&
 \\
 &&
 \exists k\forall i<k\forall j<2[n^{j}_{i}=m^{j}_{i} \,\&\, 
 (n^{1}_{k}, n^{0}_{k})<_{L}(m^{1}_{k},m^{0}_{k})]
 \\
 \mbox{ {\rm or}} && 
 \\
 &&
 \ell_{0}<\ell_{1} \,\&\, \forall i<\ell_{0}\forall j<2[n^{j}_{i}=m^{j}_{i}]
 \end{eqnarray*}
{\rm where} $<_{L}=L(<_{1},<_{0})$ {\rm denotes the lexicographic ordering:}
\[
\langle n_{1},n_{0}\rangle<_{L}\langle m_{1},m_{0}\rangle :\Leftrightarrow n_{1}<_{1}m_{1} \,\lor\,
( n_{1}=m_{1} \,\land\, n_{0}<_{0}m_{0})
.\]

{\rm Write} 
$\sum_{i<\ell}\Lam^{n^{1}_{i}}n^{0}_{i}$
{\rm for}
$\langle (n^{1}_{i}, n^{0}_{i}): i<\ell\rangle$.

\item
{\rm Let} $dom(<_{E})$ {\rm denote the domain of the relation} $<_{E}${\rm :}
\[
dom(<_{E}):=\{\sum_{i<\ell}\Lam^{n^{1}_{i}}n^{0}_{i} : \forall i<\ell\dot{-}1(n^{1}_{i+1}<_{1}n^{1}_{i}) \,\&\,
n^{1}_{i}, n^{0}_{i},\ell\in\omega\}
.\]
\item
 $<_{EW}$ {\rm denotes the restriction of} $<_{E}$ 
 {\rm to the wellfounded part in the second components:}
\begin{eqnarray*}
&&
\alpha=\sum_{i<\ell_{0}}\Lam^{n^{1}_{i}}n^{0}_{i} <_{EW} \sum_{i<\ell_{1}}\Lam^{m^{1}_{i}}m^{0}_{i} =\beta
 \mbox{ {\rm iff}}
 \\
 &&
\alpha<_{E} \beta
\,\&\, \{n^{0}_{i}:i<\ell_{0}\}\cup\{m^{0}_{i}:i<\ell_{1}\}\subseteq W(<_{0})
.\end{eqnarray*}

\end{enumerate}
\end{definition}

\bdf\label{df:LE0}
\benu

 \item\label{df:LE2}
{\rm For a definable relation} $\lhd$
{\rm and set-theoretic universe} $P$ {\rm (admissibility suffices)}
{\rm let}
\[
P\in rM_{i}(a;\lhd) :\Lrarw P\in\bigcap\{ rM_{i}(rM_{i}(b;\lhd)): b \lhd^{P} a\},
\]
{\rm where} $b\lhd^{P}a:\Lrarw P\models b\lhd a$.

{\rm Note that} $rM_{i}(a;\lhd)$ {\rm is a} $\Pi_{i+1}${\rm -class for (set-theoretic)} $\Sig_{i+1}$ $\lhd$.

\item
{\rm A relation} $\lhd$ {\rm on} $\ome$ {\rm is said to be} almost wellfounded {\rm in KP}$\ell$
{\rm if KP}$\ell$ {\rm proves the transfinite induction schema}
$TI(a,\lhd)$
{\rm up to} each $a\in\ome$.
 \eenu
\edf




\begin{lemma}\label{lem:i+1toi}{\rm (Cf. Lemma 3.2 in \cite{LMPS}.)}

Let $<_{1}, <_{0}$ be two transitive relations on $\omega$, $<_{1}$ is $\Delta_{2}$, 
$<_{0}$ is $\Sigma_{1}$, 
and $<_{EW}$ the restriction of the exponential ordering defined from these to 
the wellfounded part in the second components.
Then {\rm KP}$\ell$ proves for each $i\geq 2$
\[
\forall P\in\mbox{{\rm L}}\cup\{\mbox{{\rm L}}\}\forall a\in\omega\forall \alpha<^{P}a[
P\in rM_{i+1}(rM_{i+1}(a;<_{1}))  \to P\in rM_{i}(\alpha;<_{EW})]
\]
where for $\alpha=\sum_{i<\ell}\Lam^{n^{1}_{i}}n^{0}_{i}\in dom(<^{P}_{E})$,
$\alpha<^{P}a :\Leftrightarrow n^{1}_{0}<^{P}_{1}a$.
\end{lemma}

\begin{definition}\label{df:tower}
{\rm Let} $<_{i}\, (2\leq i\leq N-1)$ {\rm be} $\Sigma_{1}$ {\rm relations on} $\omega$.
{\rm Define a} tower {\rm relation} $<_{T}$ {\rm from these as follows.}

{\rm Define inductively relations} $<_{E_{i}}\, (2\leq i\leq N-1)$.
 \begin{enumerate}
 \item
 $<_{E_{N-1}}:\equiv <_{N-1}$.
 \item
 $<_{E_{i}}:\equiv E(<_{E_{i+1}},<_{i})$ {\rm for} $2\leq i\leq N-2${\rm , cf. Definition \ref{df:LEe}.}
 \end{enumerate}
{\rm Then let}
\[
<_{T}:\equiv <_{E_{2}}
.\]

$<_{TW}$ {\rm denotes the restriction of} $<_{T}$ {\rm to the wellfounded parts in the second components hereditarily.
Namely}
$<_{TW}=<_{E_{2}W}$ {\rm and for} $i<N-1$
\[
 \sum_{n<\ell}\Lam^{\alpha_{n}}x_{n}\in dom(<_{E_{i}W}) :\Leftrightarrow
\forall n<\ell\dot{-}1(\alpha_{n+1}<_{E_{i+1}W}\alpha_{n}) \,\&\, \forall n<\ell(x_{n}\in W(<_{i}))
\]
{\rm with} $<_{E_{N-1}W}=<_{N-1}$.

{\rm For} $a\in\omega$ {\rm and} $\alpha=\sum_{n<\ell}\Lam^{\alpha_{n}}x_{n}\in dom(<_{T})${\rm , define inductively}
\[
\alpha<a
:\Leftrightarrow
\forall n<\ell(\alpha_{n}<a)
\]
{\rm with}
$\alpha_{n}<a :\Leftrightarrow \alpha_{n}<_{N-1}a$ {\rm for}
$\alpha_{n}\in\omega$.
\end{definition}

In the following Theorem \ref{th:tower} and Corollary \ref{cor:tower},
$<_{i}\, (2\leq i\leq N-1<\omega)$ denote arbitrary
 $\Sigma_{1}$ transitive relations  on $\omega$ such that a weak theory, e.g., 
 {\rm KP}$\ell$
  proves their transitivities.

Let $<_{TW}$ denote  the restriction of the tower $<_{T}$ 
of the exponential orderings $<_{E_{i}}$
defined from these to the wellfounded parts
in the second components hereditarily.

\begin{theorem}\label{th:tower}{\rm (Cf. Theorem 3.4 in \cite{LMPS}.)}

{\rm KP}$\Pi_{N}$ proves that
\[
\forall a\in\omega\forall \alpha<a [TI(a,<_{N-1},\Pi_{N}) \to \mbox{{\rm L}}\in rM_{2}(rM_{2}(\alpha;<_{TW}))]
.\]
\end{theorem}

\bcor\label{cor:tower}
Assume that the relation $<_{N-1}$ is almost wellfounded in {\rm KP}$\ell$.
Then for {\rm each} $a\in\ome$,
\[
\mbox{{\rm KP}}\Pi_{N}\vdash
 \mbox{{\rm L}}\in \bigcap\{rM_{2}(rM_{2}(\alpha; <_{TW})): dom(<_{TW})\ni\alp<a\}
.\]
\ecor

\subsection{$k$-predecessors, relations $\prec_{k}$}\label{subsec:k-predecessor}

As in \cite{WFnonmon2} for $2\leq k\leq N$ and
ordinal terms $\alp=\psi_{\pi}^{\vec{\nu}}(a)$ with $\vec{\nu}\neq\emptyset$,
the $k$-predecessor $pd_{k}(\alp)$
is defined.
The $k$-predecessors are needed for us to embed the relation $\prec$ on $OT$ to an exponential structure
induced solely from ordinals $\{m_{k}(\alp)\}_{k}$ (cf. Lemma \ref{lem:exponentialinprec}), which in turn yields sets $V(X)=V_{N}(X)$ introduced in subsection \ref{subsec:V(X)}
with the persistency (\ref{eq:Vpersistency}) in Definition \ref{df:CX}.
As we saw it in section \ref{sec:distinguished}, the persistency is crucial for distinguished sets.

Then it turns out that $\alp\prec pd_{k}(\alp)$ holds and the $k$-predecessor $pd_{k}(\alp)$ is determined solely from
the sequences $\{\{m_{k}(\bet)\}_{2\leq k\leq N-1}:\alp\preceq\bet<\mK\}$.
Therefore it is convenient for us to handle directly the sequence of sequences $\vec{\nu}$ in defining $k$-predecessors.
After that, let us import them to ordinal terms.

Let $\pi_{i}=pd(\pi_{i+1})$ for $i<n\leq\ome$ with $\pi_{0}=\mK$.
From Definition \ref{df:notationsystem} we see that $\pi_{1}$ is defined from $\mK$ (and some $\vec{\nu},a$)
by Definition \ref{df:notationsystem}.\ref{df:notationsystem.10},
 each $\pi_{i+1}$ is defined from $\pi_{i}$ 
by Definition \ref{df:notationsystem}.\ref{df:notationsystem.11} when $1<i$ and $i\not\equiv 1 \!\!\!\!\pmod{(N-2)}$,
and  each $\pi_{i+1}$ is defined from $\pi_{i}$ 
by Definition \ref{df:notationsystem}.\ref{df:notationsystem.12} when $1<i$ and $i\equiv 1 \!\!\!\!\pmod{(N-2)}$.
In the latter case $\pi_{i+1}=\psi_{\pi_{i}}^{\vec{\nu}}(a)$ is defined from $\pi_{i},a$ and
a sequence $\vec{\nu}$ such that $\vec{\nu}<_{sp} m_{2}(\pi)$, where $m_{2}(\pi)\neq 0$ and $\fal i>2(m_{i}(\pi)=0)$.

This motivates the following.
\\

\noindent
Let $\{\vec{\nu}_{n}\}_{n\leq L}$ be a sequence of sequences $\vec{\nu}_{n}=(\nu_{n2},\ldots,\nu_{n,N-1})$ of 
ordinals
$0\neq m_{k}(\vec{\nu}_{n}):=\nu_{nk}<\veps(\Lam)$ with $0<L\equiv 0 \!\!\!\!\pmod{(N-2)}$, and $\vec{\nu}=\vec{\nu}_{0}$.

Let $pd^{(m)}(\vec{\nu}_{n})=\vec{\nu}_{n+m}$ for $n+m\leq L$.
Otherwise put $pd^{(m)}(\vec{\nu}_{n})=\emptyset$.

In what follows we assume that the following conditions are met for the sequence $\{\vec{\nu}_{n}\}_{n\leq L}$, 
$n\equiv 0 \!\!\!\!\pmod{(N-2)}$ and $2\leq k\leq N-2$.

\benu
\item
(Cf. Definition \ref{df:notationsystem}.\ref{df:notationsystem.11})
$\fal i>k(m_{i}(pd^{(k-1)}(\vec{\nu}_{n}))=0)$, 
$\fal i<k(m_{i}(pd^{(k)}(\vec{\nu}_{n}))=m_{i}(pd^{(k-1)}(\vec{\nu}_{n})))$ and
$m_{k}(pd^{(k-1)}(\vec{\nu}_{n}))=m_{k}(pd^{(k)}(\vec{\nu}_{n}))+\Lam^{m_{k+1}(pd^{(k)}(\vec{\nu}_{n}))}b$ for some 
$b<\Lam$.
In particular
\beqn\label{eq:temk}
m_{k+1}(pd^{(k)}(\vec{\nu}_{n}))=te(m_{k}(pd^{(k-1)}(\vec{\nu}_{n})))
\eeqn
and
\beqn\label{eq:N-2pd}
m_{k}(pd^{(N-2)}(\vec{\nu}_{n}))=hd(m_{k}(pd^{(k-1)}(\vec{\nu}_{n})))
\eeqn

\item
(Cf. Definition \ref{df:notationsystem}.\ref{df:notationsystem.12})
$\vec{\nu}_{n}<_{sp}m_{2}(\vec{\nu}_{n+1})$.

Let $p(\vec{\nu}_{n},m_{2}(\vec{\nu}_{n+1}))$ denote the number in Definition \ref{df:Lam}.\ref{df:Exp2.12}.

\eenu

\bdf\label{df:5resvec}
{\rm 
For $2\leq k< N$, \textit{$k$-predecessor} $pd_{k}(\vec{\nu})$ 
is defined recursively as follows.
First let $pd_{2}(\vec{\nu}):=pd(\vec{\nu})=\vec{\nu}_{1}$ with $\vec{\nu}=\vec{\nu}_{0}$.

For $k>2$, let $pd_{k}(\vec{\nu}):=\vec{\nu}_{k-1}$ if $p(\vec{\nu}_{0},m_{2}(\vec{\nu}_{1}))=0$.
Otherwise let $pd_{k}(\vec{\nu}):=pd_{k}(\vec{\nu}_{q})$ with
$q=(N-2)p(\vec{\nu}_{0},m_{2}(\vec{\nu}_{1}))$.
}
\edf

\bprp\label{prp:prdecessorvec}
$m_{k}(\vec{\nu})<_{sp}m_{k}(pd_{k}(\vec{\nu}))$.
 \eprp
\bprf
Let $\mu\leq_{pt}m_{2}(pd_{2}(\vec{\nu}))=m_{2}(pd(\vec{\nu}))$
an ordinal such that $m_{2}(\vec{\nu})<_{sd}\mu$.
$m_{2}(\vec{\nu})<_{sp}m_{2}(pd_{2}(\vec{\nu}))$ is seen from this.

Let $k>2$.
Then we have $m_{k}(\vec{\nu})<_{sd}te^{(k-2)}(\mu)$.
On the other hand we have $te^{(k-2)}(\mu)=m_{k}(pd_{k}(\vec{\nu}))$.
Hence the proposition follows.
\eprf

\bdf\label{df:5resvecnonzero}
\benu
\item\label{df:5resvecnonzero0}
{\rm Next let us define the $k$-predecessor $pd_{k}(\vec{\nu}_{i})$ for $i\not\equiv 0 \!\!\!\pmod{(N-2)}$ as follows.

Let $N-3\geq i_{0}\equiv i \!\!\!\pmod{(N-2)}$.
Then put 
$pd_{k}(\vec{\nu}_{i}):=pd(\vec{\nu}_{i})=\vec{\nu}_{i+1}$ for any $k\leq i_{0}+2$, and
$pd_{k}(\vec{\nu}_{i}):=pd^{(N-2-i_{0})}(\vec{\nu}_{i})=\vec{\nu}_{i-i_{0}+N-2}$ for $i_{0}+2<k<N$.}

\item\label{df:5resvecnonzero1}
{\rm $\vec{\nu}_{i}\prec_{k}\vec{\nu}_{j}$ denotes the transitive closure of the relation
$\{(\vec{\nu}_{i},\vec{\nu}_{j}): \vec{\nu}_{j}=pd_{k}(\vec{\nu}_{i})\}$, and
$\vec{\nu}_{i}\preceq_{k}\vec{\nu}_{j}$ its reflexive closure.}

\eenu
\edf

\bprp\label{lem:5.4vec}
Let $\vec{\mu},\vec{\xi}$ be in the sequence $\{\vec{\nu}_{n}\}_{n\leq L}$ with $\vec{\nu}_{0}=\vec{\nu}$.


Assume $\vec{\mu}\prec_{k}\vec{\xi}\prec_{k}pd_{k+1}(\vec{\mu})$. Then $pd_{k+1}(\vec{\xi})\preceq_{k}pd_{k+1}(\vec{\mu})$, and
if $pd_{k}(\vec{\mu})\neq pd_{k+1}(\vec{\mu})=pd_{k+1}(\vec{\xi})$, then $m_{k}(\vec{\mu})<_{sp}m_{k}(\vec{\xi})$.


\eprp
\bprf
Let $\vec{\mu}=\vec{\nu}_{i}$. We can assume $i\leq N-3$.

First consider the case $i\neq 0$. 
From $pd_{k}(\vec{\mu})\neq pd_{k+1}(\vec{\mu})$ 
we see that $k=i+2$,
$pd_{k+1}(\vec{\mu})=\vec{\nu}_{N-2}$, and
$pd_{k}(\vec{\mu})=\vec{\nu}_{i+1}$.
On the other side we see that $\vec{\xi}=\vec{\nu}_{N-3}$, $i=N-4$ and $k=N-2$
from $\vec{\mu}\prec_{k}\vec{\xi}$ and $pd_{k+1}(\vec{\mu})=pd_{k+1}(\vec{\xi})$.
Then $m_{N-2}(\vec{\mu})=0$ and $m_{N-2}(\vec{\mu})=0<_{sp}m_{N-2}(\vec{\xi})\neq 0$.
This shows the proposition for the case.

Next let $i=0$ and $\vec{\mu}=\vec{\nu}$. (Then $pd_{k}(\vec{\mu})\neq pd_{k+1}(\vec{\mu})$.)
Let $k>2$.
Then we have $pd_{k+1}(\vec{\mu})=pd(pd_{k}(\vec{\mu}))$ with $\vec{\xi}=pd_{k}(\vec{\mu})$.
We obtain $pd_{k+1}(\vec{\xi})=pd(\vec{\xi})=pd_{k+1}(\vec{\mu})$.
On the other hand we have $m_{k}(\vec{\mu})<_{sd}m_{k}(pd_{k}(\vec{\mu}))=m_{k}(\vec{\xi})$.
Finally let $k=2$.
We can assume $pd_{3}(\vec{\mu})>pd^{(2)}(\vec{\mu})$.
We see easily that $pd_{3}(\vec{\xi})\preceq pd_{3}(\vec{\mu})$.
If $pd_{3}(\vec{\xi})= pd_{3}(\vec{\mu})$, then we see $m_{2}(\vec{\mu})<_{sp}m_{2}(\vec{\xi})$.
\eprf

\bdf\label{df:kseqvec}
{\rm Next for $\vec{\mu}$ in the sequence $\{\vec{\nu}_{n}\}_{n\leq L}$ with $\vec{\nu}_{0}=\vec{\nu}$,
we define sequences $\{\vec{\mu}_{k}^{m}\}_{m<lh_{k}(\vec{\mu})}$ in length $lh_{k}(\vec{\mu})$ as follows.}
\benu
\item
{\rm The case when $\neg\exi\vec{\xi}(\vec{\mu}\preceq_{k}\vec{\xi}\spand pd_{k}(\vec{\xi})\neq pd_{k+1}(\vec{\xi}))$ : 
Then put} $lh_{k}(\vec{\mu})=1$ {\rm and} $\vec{\mu}_{k}^{0}:=\vec{\nu}_{L}$.

\item
{\rm The case when} $\exi\vec{\xi}(\vec{\mu}\preceq_{k}\vec{\xi}\spand pd_{k}(\vec{\xi})\neq pd_{k+1}(\vec{\xi}))${\rm : Then} 
$\vec{\mu}_{k}^{0}=\vec{\nu}_{i}$ {\rm where $i$ is the least number such that} 
$\vec{\mu}\preceq_{k}\vec{\nu}_{i}\spand pd_{k}(\vec{\nu}_{i})\neq pd_{k+1}(\vec{\nu}_{i}) ${\rm .}

{\rm Suppose that} $\vec{\mu}_{k}^{n}$ {\rm is defined so that} $pd_{k}(\vec{\mu}_{k}^{n})\neq pd_{k+1}(\vec{\mu}_{k}^{n})$.
 
 \benu
 \item 
 {\rm The case} $\exi\vec{\xi}(pd_{k+1}(\vec{\mu}_{k}^{n})\preceq_{k}\vec{\xi}\spand pd_{k}(\vec{\xi})\neq pd_{k+1}(\vec{\xi}))${\rm : Then} $\vec{\mu}_{k}^{n+1}=\vec{\nu}_{i}$ {\rm where $i$ is the least number such that} $pd_{k+1}(\vec{\mu}_{k}^{n})\preceq_{k}\vec{\nu}_{i}\spand pd_{k}(\vec{\nu}_{i})\neq pd_{k+1}(\vec{\nu}_{i})$.

 \item
 {\rm Otherwise: Then} $lh_{k}(\vec{\mu})=n+2$ {\rm and define} $\vec{\mu}_{k}^{n+1}=\vec{\nu}_{L}$.
 \eenu
\eenu
\edf

\bprp\label{lem:5Si-2vec} For $k<N-1$, 
$\vec{\mu}\preceq_{k+1}\vec{\mu}_{k}^{0}$ and 
$\fal n<lh_{k}(\vec{\mu})-1[\vec{\mu}_{k}^{n}\prec_{k+1}\vec{\mu}_{k}^{n+1}]$.
\eprp
\bprf
This is seen from the definition of $k$-predecessors in Definitions \ref{df:5resvec} and \ref{df:5resvecnonzero}.
\eprf

\bprp\label{lem:5Si-1}
Let $\vec{\nu}\prec_{k}\vec{\xi}\prec_{k}pd_{k+1}(\vec{\nu})$.
Then there exists a $\vec{\mu}\in\{\vec{\xi}\}\cup\{\vec{\xi}_{k}^{m}: m<lh_{k}(\vec{\xi})-1\}$ such that
$pd_{k+1}(\vec{\nu})=pd_{k+1}(\vec{\mu})$ and $m_{k}(\vec{\mu})>m_{k}(\vec{\nu})$.
Moreover when $\vec{\mu}\not\in\{\vec{\xi}_{k}^{m}: m<lh_{k}(\vec{\xi})-1\}$, 
$\vec{\nu}_{k}^{1}=\vec{\xi}_{k}^{0}$ holds.
\eprp
\bprf
First consider the case when $k>2$.
In the proof of Proposition \ref{lem:5.4vec} we saw
$pd_{k+1}(\vec{\nu})=pd_{k+1}(\vec{\xi})$ and $m_{k}(\vec{\xi})>m_{k}(\vec{\nu})$.
Moreover we see $\vec{\nu}_{k}^{0}=\vec{\nu}$, and $\vec{\nu}_{k}^{1}=\vec{\xi}_{k}^{0}$.

The case $k=2$ is seen from the proof of Proposition \ref{lem:5.4vec}.
\eprf

\bprp\label{lem:3.23.1vec}

Assume $\vec{\xi}=pd_{k}(\vec{\mu})$ for a $k<N-1$. Then one of the following holds:
\bdes
\item[Case \ref{lem:3.23.1vec}.1]
$\vec{\xi}=pd_{k+1}(\vec{\mu})$, $lh_{k}(\vec{\mu})=lh_{k}(\vec{\xi})$, and
\\
$\fal m<lh_{k}(\vec{\mu})[\vec{\mu}_{k}^{m}=\vec{\xi}_{k}^{m}]$.

\item[Case \ref{lem:3.23.1vec}.2]
$\vec{\mu}_{k}^{0}=\vec{\mu}$, $pd_{k+1}(\vec{\xi})=pd_{k+1}(\vec{\mu})$, $m_{k}(\vec{\xi})>m_{k}(\vec{\mu})$,
and for any 
$m<lh_{k}(\vec{\xi})=lh_{k}(\vec{\mu})-1$,
$\vec{\xi}_{k}^{m}=\vec{\mu}_{k}^{1+m}$.

\item[Case \ref{lem:3.23.1vec}.3]
$\vec{\mu}_{k}^{0}=\vec{\mu}$, $pd_{k+1}(\vec{\xi})\prec_{k}pd_{k+1}(\vec{\mu})$ and there exists an 
$m<lh(\vec{\xi})-1$ such that
$pd_{k+1}(\vec{\mu})=pd_{k+1}(\vec{\xi}^{m}_{k})$, $m_{k}(\vec{\xi}_{k}^{m})>m_{k}(\vec{\mu})$, and for any
$0<i<lh_{k}(\vec{\xi})-m=lh_{k}(\vec{\mu})$,
$ \vec{\xi}_{k}^{m+i}=\vec{\mu}_{k}^{i}$.
\edes
\eprp
\bprf
Assume $\vec{\xi}=pd_{k}(\vec{\mu})$ for a $k<N-1$. 

First consider the case $pd_{k}(\vec{\mu})= pd_{k+1}(\vec{\mu})$. 
Then $\vec{\mu}_{i}^{0}=\vec{\xi}_{i}^{0}$, and {\bf Case \ref{lem:3.23.1vec}.1} holds. 
Second suppose $pd_{k}(\vec{\mu})\neq pd_{k+1}(\vec{\mu})$. 
Then $\vec{\mu}_{k}^{0}=\vec{\mu}$ and $\vec{\mu}\prec_{k}\vec{\xi}=pd_{k}(\vec{\mu})\prec_{k}pd_{k+1}(\vec{\mu})$.
By Proposition \ref{lem:5Si-1}, if $pd_{k+1}(\vec{\xi})=pd_{k+1}(\vec{\mu})$, then {\bf Case \ref{lem:3.23.1vec}.2} holds,
i.e., {\bf Case \ref{lem:3.23.1vec}.2} holds.
Otherwise we have $pd_{k+1}(\vec{\mu})=pd_{k+1}(\vec{\xi}^{m}_{k})$ and
 $m_{k}(\vec{\xi}_{k}^{m})>m_{k}(\vec{\mu})$ for an $m<lh(\vec{\xi})-1$.
Consequently {\bf Case \ref{lem:3.23.1vec}.3} holds. 
\eprf
\\

Now let us define recursively the $k$-predecessor $pd_{k}(\alp)$ of ordinal terms $\alp=\psi_{\pi}^{\vec{\nu}}(a)$ with 
$\vec{\nu}\neq\emptyset$.

\bdf\label{df:5res}
\benu
\item
{\rm The case when $\alp=\psi_{\pi}^{\vec{\nu}}(a)$ is defined in Definition \ref{df:notationsystem}.\ref{df:notationsystem.10}.

Then $\pi=\mK$ and $\vec{\nu}=\vec{0}*(\nu)$.
Put 
$pd_{k}(\alp):=\mK$ for any $k$.}

\item
{\rm The case when $\alp=\psi_{\pi}^{\vec{\nu}}(a)$ is defined in Definition \ref{df:notationsystem}.\ref{df:notationsystem.11}.

Let $k\leq N-2$ be the number such that $ \nu_{k}=m_{k}(\pi)+\Lam^{m_{k+1}(\pi)}b$.
Then put 
$pd_{i}(\alp):=\pi$ for any $i\leq k+1$, and
$pd_{i}(\alp):=pd^{(N-k)}(\alp)$ for $k+1<i<N$, cf. Definition \ref{df:5resvecnonzero}.\ref{df:5resvecnonzero0}.
Also $pd_{N}(\alp)=\mK$.}

\item
{\rm 
The case when $\alp=\psi_{\pi}^{\vec{\nu}}(a)$ is defined in Definition \ref{df:notationsystem}.\ref{df:notationsystem.12}.

Then put $pd_{N}(\alp)=\mK$, 
$pd_{2}(\alp)=\pi$, and for $2< k\leq N-1$,
$pd_{k}(\alp):=pd^{(k-1)}(\alp)$ if $p(\vec{\nu},m_{2}(\pi))=0$.
Otherwise let $pd_{k}(\alp):=pd_{k}(pd^{(q)}(\alp))$ with
$q=(N-2)p(\vec{\nu},m_{2}(\pi))$.
}

\item
{\rm $\alp\prec_{k}\bet$ denotes the transitive closure of the relation} $\{(\alp,\bet): \bet=pd_{k}(\alp)\}$.

\eenu
\edf

\bprp\label{prp:keyfact}
Let $\sig=pd_{k+1}(\alp)\neq pd_{k}(\alp)$ for $\alp=\psi_{\pi}^{\vec{\nu}}(b)$, and 
$\pi\preceq\psi_{\sig}^{\vec{\xi}}(a)$.
Then 
$K_{\alp}(m_{k}(\alp))< a$.
\eprp
\bprf
We have $m_{k}(\alp)<_{sp}\xi_{k}$ by
Propositions \ref{prp:prdecessorvec} and \ref{lem:5.4vec}, and $\max K(m_{k}(\alp))\leq \max K(\xi_{k})$.
By (\ref{eq:notationsystem.12}) in Definition \ref{df:notationsystem}.\ref{df:notationsystem.12}, Proposition \ref{prp:G6}.\ref{prp:G6.5},
we obtain $K_{\alp}(m_{k}(\alp))<\max K(m_{k}(\alp))\leq \max K(\xi_{k})\leq a$.
\eprf

\blem\label{lem:keyfacts}
Let $\sig=pd_{k+1}(\alp)\neq pd_{k}(\alp)$ and $\bet\preceq\alp=\psi_{\pi}^{\vec{\nu}}(b)$. 

Then
$F_{\sig}(m_{k}(\alp))<\bet$.
\elem
\bprf
Let $\alp=\psi_{\pi}^{\vec{\nu}}(b)$.
By Proposition \ref{lem:5uv.1-1} we have $F_{\pi}(m_{k}(\alp))<\bet$, and it suffices to show that
$F_{\sig}(m_{k}(\alp))<\pi$ by Proposition \ref{prp:G}.\ref{prp:GF}. 
We can assume $\pi<\sig$.
Let $\pi\preceq\psi_{\sig}^{\vec{\xi}}(a)=pd_{k}(\alp)$.
By Proposition \ref{prp:keyfact} we have 
$K_{\alp}(m_{k}(\alp))< a$.
Therefore $K(m_{k}(\alp))\subset\calh_{a}(\alp)$, and
$F_{\sig}(m_{k}(\alp))\subset\calh_{a}(\pi)\cap\sig$.
Let $pd^{(i-1)}(\pi)=\pi_{i-1}=\psi_{\pi_{i}}^{\vec{\nu}_{i}}(a_{i})$
with $\pi=\pi_{0}$ and $\sig=\pi_{n}$.
We have $\calh_{a_{j+1}}(\pi_{j})\cap\pi_{j+1}\subset\pi_{j}$ and $a_{j-1}>a_{j}$ with $a=a_{n}$.
We see by induction on $n-j\geq 0$ that
$F_{\sig}(m_{k}(\alp))<\pi_{j}$.
Hence $F_{\sig}(m_{k}(\alp))<\pi_{0}=\pi$.
\eprf

\bdf\label{df:kseq}
{\rm Next for terms $\alp=\psi_{\pi}^{\vec{\nu}}(a)$ 
we define sequences $\{\alp_{k}^{m}\}_{m<lh_{k}(\alp)}$ in length $lh_{k}(\alp)$ by referring Definition \ref{df:kseqvec} as follows.}
\benu
\item
{\rm The case when $\neg\exi\del(\alp\preceq_{k}\del\spand pd_{k}(\del)\neq pd_{k+1}(\del))${\rm : Then put} $lh_{k}(\alp)=1$ {\rm and} $\alp_{k}^{0}$ {\rm is defined to be the maximal term such that} $\alp\preceq_{k+1}\alp_{k}^{0}$ {\rm with} $pd(\alp_{k}^{0})=\Lam$.}

\item
{\rm The case when} $\exi\del(\alp\preceq_{k}\del\spand pd_{k}(\del)\neq pd_{k+1}(\del))${\rm : Then} $\alp_{k}^{0}$ {\rm is defined to be the minimal term such that} $\alp\preceq_{k}\alp_{k}^{0}\spand pd_{k}(\del)\neq pd_{k+1}(\del)${\rm .}

{\rm Suppose that} $\alp_{k}^{n}$ {\rm is defined so that} $pd_{k}(\alp_{k}^{n})\neq pd_{k+1}(\alp_{k}^{n})$.
 
 \benu
 \item 
 {\rm The case} $\exi\gam(pd_{k+1}(\alp_{k}^{n})\preceq_{k}\gam\spand pd_{k}(\gam)\neq pd_{k+1}(\gam))${\rm : Then} $\alp_{k}^{n+1}$ {\rm is defined to be the minimal term such that} 
 $pd_{k+1}(\alp_{k}^{n})\preceq_{k}\alp_{k}^{n+1}$ {\rm and} 
 $pd_{k}(\alp_{k+1}^{n})\neq pd_{k+1}(\alp_{k+1}^{n})$.

 \item
 {\rm Otherwise:} $lh_{k}(\alp)=n+2$ {\rm and define} $\alp_{k}^{n+1}$ {\rm to be the maximal term such that} $\alp_{k}^{n}\preceq_{k+1}\alp_{k}^{n+1}$ {\rm with} $pd(\alp_{k}^{0})=\mK$.
 \eenu
\eenu
\edf

From Propositions \ref{lem:5Si-2vec} and \ref{lem:3.23.1vec}
we see the following Lemmas \ref{lem:5Si-2} and \ref{lem:3.23.1}.

\blem\label{lem:5Si-2} For $i<N-1$, 
$\alp\preceq_{k+1}\alp_{k}^{0}$ and 
$\fal n<lh_{k}(\alp)-1[\alp_{k}^{n}\prec_{i+1}\alp_{k}^{n+1}]$.
\elem

\blem\label{lem:3.23.1}
Assume $\eta=pd_{k}(\gam)$ for a $k<N-1$. Then one of the following holds:
\bdes
\item[Case \ref{lem:3.23.1}.1]
$\eta=pd_{k}(\gam)=pd_{k+1}(\gam)$, $lh_{k}(\gam)=lh_{k}(\eta)$, and 
$\fal m<lh_{k}(\gam)[\gam_{k}^{m}=\eta_{k}^{m}]$.

\item[Case \ref{lem:3.23.1}.2]
$\gam_{k}^{0}=\gam$, $pd_{k+1}(\eta)=pd_{k+1}(\gam)$,  $m_{k}(\eta)>m_{k}(\gam)$, 
and for any 
$m<lh_{k}(\eta)=lh_{k}(\gam)-1$,
$\eta_{k}^{m}=\gam_{k}^{1+m}$.

\item[Case \ref{lem:3.23.1}.3]
$\gam_{k}^{0}=\gam$, $pd_{k+1}(\eta)\prec_{k}pd_{k+1}(\gam)$ and there exists an 
$m<lh(\eta)-1$ such that
$pd_{k+1}(\gam)=pd_{k+1}(\eta^{m}_{k})$, $m_{k}(\eta_{k}^{m})>m_{k}(\gam_{k}^{0})$, and for any
$0<i<lh_{k}(\eta)-m=lh_{k}(\gam)$,
$ \eta_{k}^{m+i}=\gam_{k}^{i}$.

\edes
\elem

\subsection{Towers derived from ordinal terms}\label{subsec:toweronod}

In this subsection we introduce towers $T(\eta)$ of ordinal terms
 from the sequence $\{\eta_{i}^{m}:m<lh_{i}(\eta)\}$
defined in Definition \ref{df:kseq}.
We will see that the relation $\prec_{i}$ is embedded in an exponential relation $<_{E_{i}}$,
cf. Lemma \ref{lem:exponentialinprec}.

\bdf\label{df:towerordinaldiagram}

\benu
\item\label{df:towerordinaldiagram.1}
{\rm Define relations} $<_{i}$ {\rm for} $2\leq i\leq N-1$ {\rm by}
\[
\eta<_{i}\rho :\Leftrightarrow  
\eta\prec_{i}\rho \,\&\, pd_{i}(\eta)\neq pd_{i+1}(\eta)=pd_{i+1}(\rho)
\]


\item\label{df:towerordinaldiagram.4}
{\rm Extend} $<_{i}$ {\rm to} $<_{i}^{+}$ {\rm by adding the successor function} $+1${\rm .
Namely the domain is expanded to}
 $dom(<_{i}^{+}):=dom(<_{i})\cup\{a+1: a\in dom(<_{i})\}$,
{\rm and define for} $a,b\in dom(<_{i})$,
$a+1<_{i}^{+}b+1  :\Leftrightarrow  a<_{i}b$, $a+1<_{i}^{+}b  :\Leftrightarrow  a<_{i}b$, {\rm and}
$a<_{i}^{+}b+1  :\Leftrightarrow  a<_{i}b \mbox{ {\rm or} } a=b$.

$\Lam^{\alpha}$ {\rm denotes} $\Lam^{\alpha}\cdot 1$.

\item\label{df:towerordinaldiagram.3}
{\rm Let} $<_{E_{i}}$ {\rm be the exponential relation defined from}
 $<_{i}^{+}\, (2\leq i\leq N-1)$.
{\rm Namely} $<_{E_{N-1}}:\equiv <^{+}_{N-1}$ {\rm and} $<_{E_{i}}:\equiv E(<_{E_{i+1}},<^{+}_{i})$
{\rm , cf. Definition \ref{df:LEe}.}

\item\label{df:towerordinaldiagram.5}
{\rm From the sequence} $\{\eta_{i}^{m}: 2\leq i<N-1, m<lh_{i}(\eta)\}$ 
{\rm we define a tower} $T(\eta)=E_{2}(\eta)$.
{\rm The elements of the form} $E_{i}(\eta)$ {\rm are understood to be ordered by}
 $<_{E_{i}}$.
{\rm Let} $<_{T}:\equiv <_{E_{2}}$.

\begin{eqnarray*}
E_{N-1}(\eta) & := & \eta \\
E_{i}(\eta) & := & \sum_{1\leq m< lh_{i}(\eta)}\Lam^{E_{i+1}(\eta_{i}^{m})} \eta_{i}^{m-1}
+\Lam^{E_{i+1}(\eta_{i}^{0})+1} +\Lam^{E_{i+1}(\eta)}
\end{eqnarray*}

\item\label{df:towerordinaldiagram.6}
{\rm Let}
\[
\cals:=\{\la \bet,\alp\ra: \alp\preceq \bet\}
\]
{\rm On the set of pairs} $\cals$,
\[
\la x,\alp\ra<_{i,p}\la y,\bet\ra  :\Lrarw  x<^{+}_{i}y \spand \alp\prec\bet
\]
$K(\alp)=\{\alp\}$ {\rm for} $\alp\in dom(<_{E_{N-1}})$ {\rm and}
\[
K(\sum_{n<\ell}\Lam^{\alpha_{n}}x_{n})=\{x_{n}: n<\ell\}\cup\bigcup\{K(\alp_{n}): n<\ell\}
\]
$dom(<_{E_{i},p})$ {\rm is defined recursively.} $dom(<_{E_{N-1},p})=\cals$, {\rm and}
\[
\la\sum_{n<\ell}\Lam^{\alpha_{n}}x_{n},\bet\ra\in dom(<_{E_{i},p}) :\Lrarw 
\fal\gam\in K(\sum_{n<\ell}\Lam^{\alpha_{n}}x_{n})(\bet\preceq\gam)
\]
\[
\sum_{n<\ell}\Lam^{\alpha_{n}}x_{n}\in dom(<_{E_{i}}) \spand
\fal n<\ell(\la\alp_{n},\bet\ra\in dom(<_{E_{i+1},p}) \spand \bet\preceq x_{n})
\]
\[
\la\gam,\alp\ra <_{E_{i},p}\la\eta,\bet\ra  :\Lrarw  \gam<_{E_{i}}\eta \spand \alp\prec\bet
\]
{\rm and inductively define the domains} $dom(<_{E_{i}W,p})$ {\rm by} $dom(<_{E_{N-1}W,p})=\cals$,
{\rm for} $i<N-1$ {\rm and} $\sum_{n<\ell}\Lam^{\alpha_{n}}x_{n}\in dom(<_{E_{i}})$,
$ \la\sum_{n<\ell}\Lam^{\alpha_{n}}x_{n},\bet\ra\in dom(<_{E_{i}W,p})$ {\rm iff}
\[
\fal n<\ell(\la\alp_{n},\bet\ra\in dom(<_{E_{i+1}W,p}) \spand \forall n<\ell(\la x_{n},\bet\ra\in W(<_{i,p}))
\]
{\rm where} $W(<_{i,p})$ {\rm denotes the wellfounded part of} $<_{i,p}$.

$<_{T,p}:\equiv <_{E_{2,p}}$ {\rm and} $<_{TW,p} :\equiv <_{E_{2}W,p}$.

\eenu
\edf

The sequence $\{\eta_{i}^{m}: m<lh_{i}(\eta)\}$ is defined so that the following holds.

\blem\label{lem:exponentialinprec}

Suppose $\gamma\prec_{k}\eta$.
Then $\la E_{k}(\gamma),\gam\ra<_{E_{k},p}\la E_{k}(\eta),\eta\ra$.

In particular
\[
\gamma\prec_{2}\eta \Rightarrow \la T(\gamma),\gam\ra<_{T,p} \la T(\eta),\eta\ra
\]
\elem
{\bf Proof} by induction on $N-k$.

Let $\gamma\prec_{k}\eta$. It suffices to show that $E_{k}(\gamma)<_{E_{k}} E_{k}(\eta)$.
\[
E_{k}(\eta) = \sum_{1\leq n< lh_{k}(\eta)}\Lam^{E_{k+1}(\eta_{k}^{n})} \eta_{k}^{n-1}
+\Lam^{E_{k+1}(\eta_{k}^{0})+1} +\Lam^{E_{k+1}(\eta)}
\]
We can assume $\eta=pd_{k}(\gam)$.
By Lemma \ref{lem:3.23.1} one of the following cases occurs.
\bdes
\item[Case \ref{lem:3.23.1}.1]
$\eta=pd_{k}(\gam)=pd_{k+1}(\gam)$, $lh_{k}(\gam)=lh_{k}(\eta)$, and 
$\fal n<lh_{k}(\gam)[\gam_{k}^{n}=\eta_{k}^{n}]$.
Then
\[
E_{k}(\gam) = \sum_{1\leq n< lh_{k}(\eta)}\Lam^{E_{k+1}(\eta_{k}^{n})} \eta_{k}^{n-1}
+\Lam^{E_{k+1}(\eta_{k}^{0})+1} +\Lam^{E_{k+1}(\gam)}
\]

\item[Case \ref{lem:3.23.1}.2]
$\gam_{k}^{0}=\gam$, $pd_{k+1}(\eta)=pd_{k+1}(\gam)$,  $m_{k}(\eta)>m_{k}(\gam)$, 
and for any 
$n<lh_{k}(\eta)=lh_{k}(\gam)-1$,
$\eta_{k}^{n}=\gam_{k}^{1+n}$.
\[
E_{k}(\gam) = \sum_{1\leq n< lh_{k}(\eta)}\Lam^{E_{k+1}(\eta_{k}^{n})} \eta_{k}^{n-1}
+\Lam^{E_{k+1}(\eta_{i}^{0})} \gam_{k}^{0}+\Lam^{E_{k+1}(\gam_{k}^{0})+1}+\Lam^{E_{k+1}(\gam)}
\]

\item[Case \ref{lem:3.23.1}.3]
$\gam_{k}^{0}=\gam$, $pd_{k+1}(\eta)\prec_{k}pd_{k+1}(\gam)$ and there exists an 
$m<lh(\eta)-1$ such that
$pd_{k+1}(\gam)=pd_{k+1}(\eta^{m}_{k})$, $m_{k}(\eta_{k}^{m})>m_{k}(\gam_{k}^{0})$, and for any
$0<i<lh_{k}(\eta)-m=lh_{k}(\gam)$,
$ \eta_{k}^{m+i}=\gam_{k}^{i}$.
\beqnarrs
E_{k}(\eta) & = & \sum_{2\leq n< lh_{k}(\gam)}\Lam^{E_{k+1}(\gam_{k}^{n})} \gam_{k}^{n-1}
+\Lam^{E_{k+1}(\gam_{k}^{1})} \eta_{k}^{m}+E
\\
&& (E = \sum_{m\leq n< lh_{k}(\eta)}\Lam^{E_{k+1}(\eta_{k}^{n})} \eta_{k}^{n-1}
+\Lam^{E_{k+1}(\eta_{k}^{0})+1} +\Lam^{E_{k+1}(\eta)})
\\
E_{k}(\gam) & = & \sum_{2\leq n< lh_{k}(\gam)}\Lam^{E_{k+1}(\gam_{k}^{n})} \gam_{k}^{n-1}+\Lam^{E_{k+1}(\gam_{k}^{1})} \gam_{k}^{0}
+\Lam^{E_{k+1}(\gam_{k}^{0})+1} +\Lam^{E_{k+1}(\gam)}
\eeqnarrs

\edes


\eprf

\section{Wellfoundedness proof}\label{sect:wfproof}
In this section we prove Theorem \ref{th:wf}, i.e., the wellfoundedness of each initial segment of $OT$.

\subsection{The sets $V_{N}(X)$}\label{subsec:V(X)}

In this subsection sets $V(X)=V_{N}(X)$ are defined.

\bdf\label{df:5wfuv}
\benu

\item\label{df:5wfuv.Vsi}
{\rm For} $2\leq i\leq N-1$, 
\[
\bet\in U_{i}(X) :\Lrarw
[ pd_{i}(\bet)\neq pd_{i+1}(\bet)
\Rarw
F_{pd_{i+1}(\bet)}(m_{i}(\bet))\subset X]
.\]
{\rm And}
\[
\la\alp,\alp_{1}\ra<^{X}_{i,p}\la\bet,\bet_{1}\ra :\Lrarw \alp,\bet\in U_{i}(X) \spand \la\alp,\alp_{1}\ra<_{i,p}\la\bet,\bet_{1}\ra
\]
{\rm for the relation} $<_{i}$ {\rm defined in Definition \ref{df:towerordinaldiagram}.\ref{df:towerordinaldiagram.1}.}
{\rm The} domain {\rm of} $<^{X}_{i}$ {\rm is defined to be} $U_{i}(X)$.

\item
{\rm For} $2\leq i<N-1${\rm , a finite set} $\cals_{i}(\eta)$
{\rm of subterms of} $\eta$ {\rm is defines as follows:}
 \benu
 \item
 $\cals_{2}(\eta):=\{\eta_{2}^{m}: m<lh_{2}(\eta)\}$.
 \item
 {\rm For} $i>2$,
 $\cals_{i}(\eta):=\{\rho_{i}^{m}: m<lh_{i}(\rho), \rho\in\cals_{i-1}(\eta)\}$.
 \eenu

{\rm Also put} $\cals_{i}(\eta)=\emptyset$ {\rm if} $\eta$ {\rm is not of the form} $\psi_{\pi}^{\vec{\nu}}(a)$.

\item 
$\eta\in V_{N}(X)$ {\rm designates that each finite set} $\cals_{i}(\eta)$ 
{\rm  is included in the wellfounded parts} $W(<^{X}_{i})$ {\rm of the relations} $<^{X}_{i}$.
\[
\eta\in V_{N}(X)  :\Lrarw 
\fal i\in[2,N-1)\fal\bet\in\cals_{i}(\eta)[\bet\in U_{i}(X)\spand \la\bet,\eta\ra\in W(<^{X\cap\eta}_{i,p})].
\]
\eenu
\edf

It is clear that $(\bigcup\cals_{i}(\eta))\times\{\eta\}\subset\cals$ for any $\eta$, and $V_{N}(X)$ is $\Del_{1}$.
Suppose $X\cap\alp_{1}=Y\cap\alp_{1}$ and $\bet\in\cals_{i}(\eta)$ for $\eta\leq\alp_{1}$.
Then $\eta\preceq\bet$ and $F_{pd_{i+1}(\bet)}(m_{i}(\bet))<\eta$ by Lemma \ref{lem:keyfacts}.
Hence $\bet\in U_{i}(X)$ iff $\bet\in U_{i}(Y)$.
Obviously $\la\alp,\gam\ra<_{i}^{X\cap\eta}\la\bet,\eta\ra \Lrarw \la\alp,\gam\ra<_{i}^{Y\cap\eta}\la\bet,\eta\ra$
since $F_{pd_{i+1}(\alp)}(m_{i}(\alp))<\gam\leq\eta$ by Lemma \ref{lem:keyfacts}
and 
$\gam\preceq\alp$, $\gam\prec\eta$.
Therefore $\la\bet,\eta\ra\in W(<_{i,p}^{X\cap\eta})$ iff $\la\bet,\eta\ra\in W(<_{i,p}^{Y\cap\eta})$.
Thus $V_{N}(X)$ enjoys the condition (\ref{eq:Vpersistency}).

\bprp\label{lem:gvinclu}
For any limit universe $P$,
if $\gam\in\calg(\calw^{P})$, then $\fal i\in[2,N-1)[\cals_{i}(\gam)\subset U_{i}(\calw^{P})]$
and $\cals_{N-2}(\gam)\subset U_{N-1}(\calw^{P})$.
\eprp
\bprf
Assume $\gam\in\calg(\calw^{P})$.
Let $\del\in\cals_{i}(\gam)$, $\nu=m_{i}(\del)$ and $\sig=pd_{i+1}(\del)$.
Then $\gam\preceq\del$.
We have to show $F_{\sig}(\nu)\subset \calw^{P}$.
By Lemma \ref{lem:keyfacts} 
we have $F_{\sig}(\nu)<\gam$.

On the other hand we have $\gam\in\calc^{\gam}(\calw^{P})$, and this yields 
 $\nu\in\calc^{\gam}(\calw^{P})$ by the definition of the set $\calc^{\gam}(\calw^{P})$.
Therefore $F_{\sig}(\nu)\subset\calc^{\gam}(\calw^{P})$ follows from Proposition \ref{lem:CMsmpl}.
Thus we have $F_{\sig}(\nu)\subset\calc^{\gam}(\calw^{P})\cap\gam\subset \calw^{P}$.

For the case $i=N-2$, let $\mu=m_{N-1}(\del)$ with $\mK=pd_{N}(\del)$ and $\del\in\cals_{N-2}(\gam)$.
$F_{\mK}(\mu)\subset\calw^{P}\cap\gam$ is seen from $F_{\mK}(\mu)<\gam$.
\eprf
\\

Let $<^{P}_{T\calw,p}$ denote $(<_{T\calw,p})^{P}$, i.e., the relation $<_{T\calw,p}=<_{E_{2}\calw,p}$ in $P$
\blem\label{lem:TWelementary}
\benu

\item\label{lem:TWelementary3}
$<_{N-1}$ is almost wellfounded in {\rm KP}$\ell$.



\item\label{lem:TWelementary5}
Let $P$ be a limit universe.
Suppose $\eta\in V_{N}(\calw^{P})$.
Then $\la T(\eta),\eta\ra\in dom(<^{P}_{T\calw,p})$.
Moreover if $\gam\prec\eta$ and $\gam\in V_{N}(\calw^{P})$,
then $\la T(\gam),\gam\ra<^{P}_{T\calw,p}\la T(\eta),\eta\ra$.

\eenu
\elem
\bprf
\\
\ref{lem:TWelementary}.\ref{lem:TWelementary3}.
$\gam<_{N-1}\eta\Lrarw \gam\prec_{N-1}\eta$, and this implies
$m_{N-1}(\gam)<m_{N-1}(\eta)<\Lam=\veps_{\mK+1}$.
\\

\noindent
\ref{lem:TWelementary}.\ref{lem:TWelementary5}. 
The fact that $\eta\in V_{N}(\calw^{P})\Rarw \la T(\eta),\eta\ra\in dom(<^{P}_{T\calw,p})$ is seen from the definition of 
$<^{P}_{T\calw,p}$.
Assume $\gam\prec\eta$ and $\gam\in V_{N}(\calw^{P})$.
Then by Lemma \ref{lem:exponentialinprec} we have $\la T(\gam),\gam\ra<_{T,p}\la T(\eta),\eta\ra$.
Moreover we have $\la T(\gam),\gam\ra\in dom(<^{P}_{T\calw,p})$.
Hence $\la T(\gam),\gam\ra<^{P}_{T\calw,p}\la T(\eta),\eta\ra$.
\eprf

\blem\label{lem:3wf16}

If $P\in rM_{2}(rM_{2}(T(\eta) ;<_{TW,p}))$, then
$\eta\in\calg(\calw^{P})\cap V_{N}(\calw^{P}) \to \eta\in\calw^{P}$.
\elem
{\bf Proof} by induction on $\in$.


Let $\calx=rM_{2}(T(\eta) ;<_{TW,p})\subset Lmtad$.
First we show the existence of a distinguished set $X_{1}\in P$ such that
\beqn\renewcommand{\theequation}{\ref{eq:Vlocalize}}
\fal Q\in P\cap\calx[X_{1}\in Q \Rarw \eta\in V_{N}(\calw^{Q})]
\eeqn
\addtocounter{equation}{-1}
We have $\fal i\in[2,N-1)\fal\bet\in\cals_{i}(\eta)[F_{pd_{i+1}(\bet)}(m_{i}(\bet))\subset\calw^{P}\cap\eta]$.
Pick a distinguished set $X_{1}\in P$ such that 
$\fal i\in[2,N-1)\fal\bet\in\cals_{i}(\eta)[F_{pd_{i+1}(\bet)}(m_{i}(\bet))\subset X_{1}\cap\eta]$.
Let $X_{1}\in Q\in P\cap\calx$.
Then $X_{1}\subset\calw^{Q}\subset\calw^{P}$, and hence 
$\fal i\in[2,N-1)\fal\bet\in\cals_{i}(\eta)[F_{pd_{i+1}(\bet)}(m_{i}(\bet))\subset\calw^{Q}\cap\eta]$, i.e.,
$\fal i\in[2,N-1)\fal\bet\in\cals_{i}(\eta)[\bet\in U_{i}(\calw^{Q}\cap\eta)]$

Furthermore we have $\bet\in W(<_{i}^{\calw^{P}\cap\eta})$ for $\bet\in\cals_{i}(\eta)$,
  and $\calw^{Q}\subset\calw^{P}$.
Hence $U_{i}(\calw^{Q}\cap\eta)\subset U_{i}(\calw^{P}\cap\eta)$ and $\bet\in W(<_{i}^{\calw^{Q}\cap\eta})$.
We obtain $\eta\in V_{N}(\calw^{Q})$.

By Corollary \ref{cor:3wf16} it suffices to show (\ref{eq:etamahlo}) for any $Q\in P\cap\calx$
such that $X_{1}\in Q$.

\beqn\renewcommand{\theequation}{\ref{eq:etamahlo}}
\forall\gamma\prec\eta\{\gam\in\calg(\calw^{Q})\cap V_{N}(\calw^{Q})
 \Rightarrow   \gamma\in \calw^{Q}\}
\end{equation}
\addtocounter{equation}{-1}

Let $Q\in P\cap\calx$, $X_{1}\in Q$ and assume
$\gam\prec\eta$ and $\gam\in\calg(\calw^{Q})\cap V_{N}(\calw^{Q})$.
Then $T(\gam)<^{Q}_{TW;\gam}T(\eta)$
by Lemma \ref{lem:TWelementary}.\ref{lem:TWelementary5}.

Therefore $Q\in rM_{2}(rM_{2}(T(\gam) ;<_{TW,p}))$
by $Q\in\calx=rM_{2}(T(\eta) ;<_{TW,})$.
IH on $\in$ yields $\gamma\in \calw^{Q}$. 
This shows (\ref{eq:etamahlo}).
We conclude $\eta\in\calw^{P}$ by Corollary \ref{cor:3wf16}.
\eprf

\blem\label{th:GWconcl}
For {\rm each} $n\in\ome$ 
\[
{\sf KP}\Pi_{N}\vdash \fal\alp\in OT_{n}[\alp\in\calg(\calw)\cap V_{N}(\calw)\cap\mK \to \alp\in \calw]
.\]
\elem
\bprf
This is seen from Proposition \ref{lem:4acalg}, Corollary \ref{cor:tower}  and Lemmas \ref{lem:TWelementary}.\ref{lem:TWelementary3} and \ref{lem:3wf16}.
\eprf

\subsection{Wellfoundedness proof (concluded)}\label{subsec:wfproofconclude}
In the final subsection we conclude the wellfoundedness.

Let for $\vec{\xi}=(\xi_{2},\ldots,\xi_{N-1})$ with $\xi_{i}\in E$
\beqnarrs
E_{n}\calc^{\mK}(\calw) & := & \{\xi\in E : K(\xi)\subset\calc^{\mK}(\calw)\cap OT_{n}\}
\\
\vec{E}_{n}\calc^{\mK}(\calw) & := & \{\vec{\xi}\subset E_{n}\calc^{\mK}(\calw): \vec{\xi} \mbox{ {\rm is irreducible}}\}
\eeqnarrs

\bdf\label{df:id4wfA}{\rm For} $a\in OT_{n}$ {\rm and irreducible sequences} $\vec{\nu}=(\nu_{2},\ldots,\nu_{N-1})\subset E_{n}${\rm , define:}
\benu
\item\label{df:id4wfA.1}
\[
 A(a,\vec{\nu}) :\Lrarw 
 \fal\sig\in\calw\cup\{\mK\}[\psi_{\sig}^{\vec{\nu}}(a)\in OT_{n} \Rarw \psi_{\sig}^{\vec{\nu}}(a)\in\calw].
\]
\item \label{df:id4wfA.2}
\[
\mbox{{\rm MIH}}(a) :\Lrarw
 \fal b\in\calc^{\mK}(\calw)\cap a\fal \vec{\nu}\in\vec{E}_{n}\calc^{\mK}(\calw)\, A(b,\vec{\nu}).
\]
\item\label{df:id4wfA.3}
\[
\mbox{{\rm SIH}}(a,\vec{\nu}) :\Lrarw
 \fal \vec{\xi}\in\vec{E}_{n}\calc^{\mK}(\calw)[\vec{\xi}<_{lx,2}\vec{\nu} \Rarw A(a,\vec{\xi})]
.\]
\eenu
\edf

\blem\label{th:id5wf21}
Assume $\{a\}\cup K(\vec{\xi})\subset\calc^{\mK}(\calw)$, $\mbox{{\rm MIH}}(a)$, and $\mbox{{\rm SIH}}(a,\vec{\xi})$ in Definition \ref{df:id4wfA}.
Then
\[
 \fal\kap\in\calw\cup\{\mK\}[\psi_{\kap}^{\vec{\xi}}(a)\in OT_{n} \Rarw \psi_{\kap}^{\vec{\xi}}(a)\in\calg(\calw)].
\]
\elem
\bprf  
Let $\alp_{1}=\psi_{\kap}^{\vec{\xi}}(a)\in OT_{n}$ with $\kap\in\calw\cup\{\mK\}$. We have to show $\alp_{1}\in\calg(\calw)$.

By Proposition \ref{lem:CX2}.\ref{lem:CX2.3} we have 
$\{\kap,a\}\cup K(\vec{\xi})\subset\calc^{\mK}(\calw)\subset\calc^{\alp_{1}}(\calw)$
 and hence by Lemma \ref{lem:wf5.332}
\beqn\label{eqn:id4wf21.1}
\alp_{1}\in\calc^{\alp_{1}}(\calw) \spand \fal\rho[G_{\rho}(\{\kap,a\}\cup K(\vec{\xi}))\subset\calw]
\eeqn

Thus it suffices to show the following claim.
\bclm\label{clm:id5wf21.1}
$$\fal\bet_{1}\in\calc^{\alp_{1}}(\calw)\cap\alp_{1}[\bet_1\in\calw].$$
\eclm
{\bf Proof} of Claim \ref{clm:id5wf21.1} by induction on $\ell\bet_1$. 
Assume $\bet_{1}\in\calc^{\alp_{1}}(\calw)\cap\alp_{1}$ and let
\[
\mbox{LIH} :\Lrarw
\fal\gam\in\calc^{\alp_{1}}(\calw)\cap\alp_{1}[\ell\gam<\ell\bet_{1} \Rarw \gam\in\calw].
\]

We show $\bet_1\in\calw$. 


\noindent
{\bf Case 0}. 
$\bet_1\not\in\calE(\bet_{1})$ or $\bet_{1}\in\calw\cap\alp_{1}$: 
Assume $\bet_{1}\not\in\calw$.
Then $S(\bet_{1})\subset\calc^{\alp_{1}}(\calw)\cap\alp_{1}$. 
LIH yields $S(\bet_{1})\subset\calw$. 
Hence we conclude $\bet_{1}\in\calw$ from Proposition \ref{lem:id3wf19-1}.
\\

In what follows consider the cases when $\bet_{1}=\psi_{\pi}^{\vec{\nu}}(b)$ for some $\pi,b,\vec{\nu}$.
We can assume $\{\pi,b\}\cup K(\vec{\nu})\subset\calc^{\alp_{1}}(\calw)$.
\\
{\bf Case 1}. $\pi\leq\alp_{1}$: 
Then $\{\bet_{1}\}=G_{\pi}(\bet_{1})\subset\calw$ by $\bet_{1}\in\calc^{\alp_{1}}(\calw)$ and Proposition \ref{lem:KC}.
%
\\

\noindent
{\bf Case 2}.
$b<a$, $\bet_{1}<\kap$ and $K_{\alp_{1}}(\{\pi,b\}\cup K(\vec{\nu}))<a$:
Let $B$ denote a set of subterms of $\bet_{1}$ defined recursively as follows.
First $\{\pi,b\}\cup K(\vec{\nu})\subset B$.
Let $\alp_{1}\leq\bet\in B$. 
If $\bet=_{NF}\ome^{\gam}>\mK$, then $\gam\in B$.
If $\bet=_{NF}\gam_{m}+\cdots+\gam_{0}$, then $\{\gam_{i}:i\leq m\}\subset B$.
If $\bet=_{NF}\vphi\gam\del$, then $\{\gam,\del\}\subset B$.
If $\bet=_{NF}\Ome_{\gam}$, then $\gam\in B$.
If $\bet=_{NF}\psi_{\sig}^{\vec{\zeta}}(c)$, then $\{\sig,c\}\cup K(\vec{\zeta})\subset B$.

Then from $\{\pi,b\}\cup K(\vec{\nu})\subset\calc^{\alp_{1}}(\calw)$ we see inductively that
$B\subset\calc^{\alp_{1}}(\calw)$.
Hence by LIH we have $B\cap\alp_{1}\subset\calw$.
Moreover if $\alp_{1}\leq\psi_{\sig}^{\vec{\zeta}}(c)\in B$, then $c\in K_{\alp_{1}}(\{\pi,b\}\cup K(\vec{\nu}))<a$.

We claim that
\bclm\label{eq:case2A}
$\fal\bet\in B(\bet\in\calc^{\mK}(\calw))$.
\eclm
{\bf Proof} of Claim \ref{eq:case2A} by induction on $\ell\bet$.
Let $\bet\in B$. We can assume that $\alp_{1}\leq\bet=\psi_{\sig}^{\vec{\zeta}}(c)$ by induction hypothesis on the lengths.
Then by induction hypothesis we have $\{\sig,c\}\cup K(\vec{\zeta})\subset\calc^{\mK}(\calw)$.
On the other hand we have $c<a$.
$\mbox{MIH}(a)$ yields $\bet\in\calw$.
Thus the Claim \ref{eq:case2A} is shown.
\eprf
\\

In particular we obtain $\{\pi,b\}\cup K(\vec{\nu})\subset \calc^{\mK}(\calw)$.
Moreover we have $b<a$.
Therefore once again $\mbox{MIH}(a)$ yields $\bet_{1}\in\calw$.
\\

\noindent
{\bf Case 3}. 
$b=a$, $\pi=\kap$, $\fal\del\in K(\vec{\nu})(K_{\alp_{1}}(\del)<a)$ and $\vec{\nu}<_{lx,2}\vec{\xi}$: 
As in Claim \ref{eq:case2A} we see that $K(\vec{\nu})\subset\calc^{\mK}(\calw)$ from $\mbox{MIH}(a)$.
$\mbox{SIH}(a,\vec{\xi})$ yields $\bet_{1}\in\calw$.
\\

\noindent
{\bf Case 4}.
$a\leq b\leq K_{\bet_{1}}(\del)$ for some $\del\in K(\vec{\xi})\cup\{\kap,a\}$:
It suffices to find a $\gam$ such that $\bet_{1}\leq\gam\in\calw\cap\alp_{1}$.
Then $\bet_{1}\in\calw$ follows from $\bet_{1}\in\calc^{\alp_{1}}(\calw)$ and Proposition \ref{prp:updis}.

We see that $a\in K_{\del}(\alp)$ iff $\psi_{\kap}^{\vec{\xi}}(a)\in k_{\del}(\alp)$ for some $\kap,\vec{\xi}$,
and for each $\psi_{\kap}^{\vec{\xi}}(a)\in k_{\del}(\psi_{\kap_{0}}^{\vec{\xi}_{0}}(a_{0}))$ there exists a sequence
$\{\alp_{i}\}_{i\leq m}$ of subterms of $\alp_{0}=\psi_{\kap_{0}}^{\vec{\xi}_{0}}(a_{0})$ such that 
$\alp_{m}=\psi_{\kap}^{\vec{\xi}}(a)$, 
$\alp_{i}=\psi_{\kap_{i}}^{\vec{\xi}_{i}}(a_{i})$ for some $\kap_{i},a_{i},\vec{\xi}_{i}$,
and for each $i<m$,
$\del\leq\alp_{i+1}\in\calE(C_{i})$ for $C_{i}=\{\kap_{i},a_{i}\}\cup K(\vec{\xi}_{i})$.

Pick an $\alp_{2}=\psi_{\kap_{2}}^{\vec{\xi}_{2}}(a_{2})\in \calE(\del)$ and 
an $\alp_{m}=\psi_{\kap_{m}}^{\vec{\xi}_{m}}(a_{m})\in k_{\bet_{1}}(\alp_{2})$ for some $\kap_{m},\vec{\xi}_{m}$ and
$a_{m}\geq b\geq a$.
We have $\alp_{2}\in\calw$ by $\del\in\calc^{\mK}(\calw)$.
We can assume $\alp_{2}\geq\alp_{1}$.
Then $a_{2}\in K_{\alp_{1}}(\alp_{2})<a\leq b$, and $m>2$.

Let $\{\alp_{i}\}_{2\leq i\leq m}$ be the sequence of subterms of $\alp_{2}$ such that
$\alp_{i}=\psi_{\kap_{i}}^{\vec{\xi}_{i}}(a_{i})$ for some $\kap_{i},a_{i},\vec{\xi}_{i}$,
and for each $i<m$,
$\bet_{1}\leq\alp_{i+1}\in\calE(C_{i})$ for $C_{i}=\{\kap_{i},a_{i}\}\cup K(\vec{\xi}_{i})$.
Let $\{n_{j}\}_{0\leq j\leq k}\, (0<k\leq m-2)$ be the increasing sequence $n_{0}<n_{1}<\cdots<n_{k}\leq m$ 
defined recursively by $n_{0}=2$, and assuming $n_{j}$ has been defined so that
$n_{j}<m$ and $\alp_{n_{j}}\geq\alp_{1}$, $n_{j+1}$ is defined as follows
\[
n_{j+1}=\min(\{i: n_{j}\leq i<m: \alp_{i}<\alp_{n_{j}}\}\cup\{m\})
.\]
If either $n_{j}=m$ or $\alp_{n_{j}}<\alp_{1}$, then $k=j$ and $n_{j+1}$ is undefined.

Then we claim that
\bclm\label{eq:case4A}
$\fal j\leq k(\alp_{n_{j}}\in\calw) \spand \alp_{n_{k}}<\alp_{1}$.
\eclm
{\bf Proof} of Claim \ref{eq:case4A}.
By induction on $j\leq k$ we show first that $\fal j\leq k(\alp_{n_{j}}\in\calw)$. 
We have $\alp_{n_{0}}=\alp_{2}\in\calw$.
Assume $\alp_{n_{j}}\in\calw$ and $j<k$.
Then $n_{j}<m$, i.e., $\alp_{n_{j+1}}<\alp_{n_{j}}$, and 
by $\alp_{n_{j}}\in\calc^{\alp_{n_{j}}}(\calw)$, we have $C_{n_{j}}\subset\calc^{\alp_{n_{j}}}(\calw)$,
and hence $\alp_{n_{j}+1}\in\calE(C_{n_{j}})\subset\calc^{\alp_{n_{j}}}(\calw)$.
We see inductively that
$\alp_{i}\in \calc^{\alp_{n_{j}}}(\calw)$ for any $i$ with $n_{j}\leq i\leq n_{j+1}$.
Therefore $\alp_{n_{j+1}}\in \calc^{\alp_{n_{j}}}(\calw)\cap\alp_{n_{j}}\subset\calw$ by Proposition \ref{prp:maxwup}.

Next we show that $\alp_{n_{k}}<\alp_{1}$.
We can assume that $n_{k}=m$.
This means that $\fal i(n_{k-1}\leq i<m \Rarw \alp_{i}\geq\alp_{n_{k-1}})$.
We have
$\alp_{2}=\alp_{n_{0}}>\alp_{n_{1}}>\cdots>\alp_{n_{k-1}}\geq\alp_{1}$, and
$\fal i<m(\alp_{i}\geq\alp_{1})$.
Therefore $\alp_{m}\in k_{\alp_{1}}(\alp_{2})\subset k_{\alp_{1}}(\{\kap,a\}\cup K(\vec{\xi}))$, i.e.,
$a_{m}\in K_{\alp_{1}}(\{\kap,a\}\cup K(\vec{\xi}))$ for $\alp_{m}=\psi_{\kap_{m}}^{\vec{\xi}_{m}}(a_{m})$.
On the other hand we have $K_{\alp_{1}}(\{\kap,a\}\cup K(\vec{\xi}))<a$ for 
$\alp_{1}=\psi_{\kap}^{\vec{\xi}}(a)$.
Thus $a\leq a_{m}<a$, a contradiction.

The Claim \ref{eq:case4A} is shown, and we obtain $\bet_{1}\leq\alp_{n_{k}}\in\calw\cap\alp_{1}$.

This completes a proof of Claim \ref{clm:id5wf21.1} and of the lemma.
\eprf

\blem\label{lem:id3wf20}
Suppose $\mbox{{\rm MIH}}(a)$ and $\kap\leq\mK$.
 For any ordinal term $\bet\in OT_{n}$
\[
F_{\kap}(\bet)\subset\calw \spand K_{\kap}(\bet)<a
 \Rarw \bet\in\calc^{\mK}(\calw).
\]
\elem
{\bf Proof} by induction on $\ell\bet$.
By IH with Proposition \ref{lem:id3wf19-1} we can assume $\bet=\psi_{\rho}^{\vec{\nu}}(b)\geq\kap$.
Then $F_{\kap}(\bet)=F_{\kap}(\{\rho,b\}\cup K(\vec{\nu}))$ and
$\{b\}\cup K_{\kap}(\{\rho,b\}\cup K(\vec{\nu}))=K_{\kap}(\bet)<a$.
By IH we have $\{\rho,b\}\cup K(\vec{\nu})\subset\calc^{\mK}(\calw)$.
$\mbox{{\rm MIH}}(a)$ with $b<a$ yields $A(b,\vec{\nu})$, and we obtain
$\bet=\psi_{\rho}^{\vec{\nu}}(b)\in\calw$ by $\rho\in\calw\cup\{\mK\}$.
\eprf


\bprp\label{lem:TI}
For {\rm each} $n<\ome$, ${\sf KP}\ell\vdash TI[\calc^{\mK}(\calw)\cap\ome_{n+1}(\mK+1)]$.
\eprp
\bprf
 By metainduction on $n<\ome$ using Proposition \ref{lem:3wf6}
we see $TI[\calc^{\mK}(\calw)\cap\ome_{n+1}(\mK+1)]$, i.e., 
$Prg[\calc^{\mK}(\calw)\cap\ome_{n+1}(\mK+1),\caly]\to \calc^{\mK}(\calw)\cap\ome_{n+1}(\mK+1)\subset\caly$ for any definable class $\caly$.
\eprf

\blem\label{th:id5wf21V}
Assume $\{a\}\cup K(\vec{\xi})\subset\calc^{\Lam}(\calw)$, $\mbox{{\rm MIH}}(a)$, and $\mbox{{\rm SIH}}(a,\vec{\xi})$ in Definition \ref{df:id4wfA}.
Then
\[
 \fal\pi\in\calw\cup\{\mK\}[\psi_{\pi}^{\vec{\xi}}(a)\in OT_{n} \Rarw \psi_{\pi}^{\vec{\xi}}(a)\in V_{N}(\calw)].
\]
\elem
\bprf
By Lemmas \ref{th:id5wf21} and \ref{th:GWconcl} it suffices to show that
$\alp_{1}=\psi_{\pi}^{\vec{\xi}}(a)\in V_{N}(\calw)$, cf. Definition \ref{df:5wfuv}.
Let $2\leq i<N-1$, $\bet_{1}=\psi_{\sig}^{\vec{\mu}}(b)\in\cals_{i}(\alp_{1})$.
We have to show $\bet_{1}\in W(<^{\calw\cap\alp_{1}}_{i})$.
Suppose $pd_{i}(\bet_{1})\neq pd_{i+1}(\bet_{1})$ and $\gam_{1}<^{\calw\cap\alp_{1}}_{i}\bet_{1}$.
We have $\gam_{1}\in U_{i}(\calw\cap\alp_{1})$, and $\gam_{1}<_{i}\bet_{1}$, i.e., $\gam_{1}\prec_{i}\bet_{1}$,
$pd_{i}(\gam_{1})\neq pd_{i+1}(\gam_{1})$, and $\kap:=pd_{i+1}(\gam_{1})=pd_{i+1}(\bet_{1})$,
cf. Definition \ref{df:towerordinaldiagram}.\ref{df:towerordinaldiagram.1}.
Hence $\nu:=m_{i}(\gam_{1})<m_{i}(\bet_{1})$ by Proposition \ref{lem:5.4vec}.

We claim that $\nu\in\calc^{\mK}(\calw)$.
By Lemma \ref{lem:id3wf20} it suffices to show that
$F_{\kap}(\nu)\subset\calw$ and $K_{\kap}(\nu)<a$.
We have $F_{\kap}(\nu)\subset\calw$ by $\gam\in U_{i}(\calw\cap\alp_{1})$.

By Proposition \ref{prp:G6}.\ref{prp:G6.3} and $\bet_{1}\prec\kap$ we have $c\leq b\leq a$.
On the other hand we have 
$K_{\gam_{1}}(\nu)<c$ for $c=a_{i}(\gam_{1})$, i.e., $K(\nu)\subset\calh_{c}(\gam_{1})$
by Proposition \ref{prp:keyfact}.
Then $K(\nu)\subset\calh_{a}(\kap)$ by $a\geq c$ and $\kap>\gam_{1}$, and
$K_{\kap}(\nu)<a$.

Thus we have shown $\nu=m_{i}(\gam_{1})\in\calc^{\mK}(\calw)$.
Therefore $\bet_{1}\in W(<^{\calw\cap\alp_{1}}_{i})$ is seen by induction on 
$m_{i}(\gam_{1})\in\calc^{\mK}(\calw)\cap\ome_{n}(\mK+1)$,
cf. Proposition \ref{lem:TI}.
\eprf

\bprp\label{lem:WfCKWP}
For {\rm each} $\in\ome$ and each definable class $\calx$ of irreducible sequences $\vec{\xi}=(\xi_{2},\ldots,\xi_{N-1})$
of $\xi_{i}<\ome_{n}(\mK+1)$
\[
{\sf KP}\ell\vdash Prg_{lx}[\vec{E}_{n}\calc^{\mK}(\calw),\calx] \to
\fal\vec{\xi}\in\vec{E}_{n}\calc^{\mK}(\calw)(\vec{\xi}\in\calx)
\]
where
\[
Prg_{lx}[\vec{E}_{n}\calc^{\mK}(\calw),\calx] :\Lrarw
\fal\vec{\xi}\in\vec{E}_{n}\calc^{\mK}(\calw)[\fal \vec{\nu}\in\vec{E}_{n}\calc^{\mK}(\calw)(\vec{\nu}<_{lx,2}\vec{\xi} \to \vec{\nu}\in\calx) \to \vec{\xi}\in\calx]
.\]

\eprp
\bprf
In Definition \ref{df:lxo} ordinals $o(\vec{\xi})<\veps(\Lam)=\veps_{\mK+2}$ are assigned to irreducible $\vec{\xi}$ so that
$\vec{\nu}<_{lx,2}\vec{\xi} \Rarw o(\vec{\nu})<o(\vec{\xi})$ by Proposition \ref{prp:lxo}, and
$K(o(\vec{\xi}))\subset\calc^{\mK}(\calw)$ if $\vec{\xi}\in\vec{E}_{n}\calc^{\mK}(\calw)$.

Now since $K(\vec{\xi})<\ome_{n}(\mK+1)$, we can replace each occurrence of $\Lam=\veps_{\mK+1}$ in $\vec{\xi}$
by $\lam_{n}:=\ome_{n}(\mK+1)$: let $o_{n}(\vec{\nu})$ denote the result of replacing $\Lam$ by $\lam_{n}$
in $o(\vec{\nu})$. Then
$\vec{\nu}<_{lx,2}\vec{\xi} \Rarw o_{n}(\vec{\nu})<o_{n}(\vec{\xi})$ for any $\vec{\nu},\vec{\xi}$ such that
$K(\{\vec{\nu},\vec{\xi}\})<\lam_{n}$.

Furthermore we have $o_{n}(\vec{\xi})<\ome_{n(N-1)}(\mK+1)$ since
$\mK\cdot\ome_{n}(\mK+1)=\ome_{n}(\mK+1)$ for $n>1$.
Hence the proposition follows from Proposition \ref{lem:TI}.
\eprf
\\

\noindent
Using Lemma \ref{th:id5wf21V}, Propositions \ref{lem:TI} and \ref{lem:WfCKWP} we see
\[
\fal a\in\calc^{\mK}(\calw)\cap\ome_{n}(\mK+1)\fal\vec{\nu}\in\vec{E}_{n}\calc^{\mK}(\calw)\, A(a,\vec{\nu})
\]
by main induction on $a\in\calc^{\mK}(\calw)\cap\ome_{n}(\mK+1)$ with subsidiary induction on 
$\vec{\xi}\in\vec{E}_{n}\calc^{\mK}(\calw)\, (K(\vec{\xi})<\ome_{n}(\mK+1))$ along $<_{lx,2}$.
Hence by induction on $\ell\alp$ we see that $\alp\in OT_{n}\Rarw\alp\in\calc^{\mK}(\calw)$.
Thus Theorem \ref{lem:wfTn}, and hence Theorem \ref{th:wf} is shown.

\end{document}